\documentclass[10pt,a4paper]{article}

\usepackage[hidelinks]{hyperref}
\usepackage{lineno}

\usepackage{amssymb, bm}
\usepackage{amsmath}
\usepackage{amsthm}
\usepackage{mathrsfs}
\usepackage{siunitx}
\usepackage[margin=2cm]{geometry}
\usepackage{subfig}
\usepackage{float}
\usepackage{xcolor}
\usepackage{subfloat}

\usepackage{mathtools}
\usepackage{bbm}
\usepackage{multirow}
\usepackage{stmaryrd}

\usepackage{algorithm}     
\usepackage{algpseudocode} 

\usepackage[utf8]{inputenc} 
\usepackage[T1]{fontenc}

\usepackage{tikz}
\usetikzlibrary{decorations.markings}
\usetikzlibrary{shapes.geometric}
\usetikzlibrary{shapes.arrows}
\usetikzlibrary{fit,calc}
\usetikzlibrary{quotes,angles}
\usetikzlibrary{arrows.meta}
\usepackage{animate}  
\usepackage{graphicx}
\usepackage{pgf}

\usepackage{pgfplots}
\usepackage[normalem]{ulem}

\newcommand{\vertiii}[1]{{\left\vert\kern-0.25ex\left\vert\kern-0.25ex\left\vert #1 
   \right\vert\kern-0.25ex\right\vert\kern-0.25ex\right\vert}}

\usepackage{geometry}
\makeatletter
\def\blfootnote{\xdef\@thefnmark{$\star$}\@footnotetext}
\makeatother
\newenvironment{Authors}%
  {\begin{center}\begin{bfseries}}%
  {\end{bfseries}\end{center}}
\newenvironment{Addresses}%
  {\begin{flushleft}\begin{itshape}}%
  {\end{itshape}\end{flushleft}}
  \newcommand{\email}[1]{\hspace*{\stretch{1}}\emph{\texttt{#1}}}

 \usepackage{fancyhdr}  
  \fancypagestyle{plain}{
\fancyhead{}
\fancyhead[C]{\hfill Submitted to Elsevier, Mars 2025}
 }

\newtheorem{definition1}{Definition}
\newtheorem{theorem}{Theorem}[section]
\newtheorem{theorem1}{Theorem}

\newtheorem{lemma1}[theorem1]{Lemma}
\newtheorem{remark}[theorem]{Remark}
\newtheorem*{remark*}{Remark}  

\begin{document}

\title{Optimization-based method for conjugate heat transfer problems}

 \date{}
 \maketitle
\vspace{-50pt} 
 
\begin{Authors}
Liang Fang$^{1}$,  Xiandong Liu$^{2}$,
Lei Zhang$^2$.
\end{Authors}

\begin{Addresses}
$^1$
School of Automotive Studies, Tongji University,
Shanghai 201804, China, \email{fangliang@tongji.edu.cn} \\[3mm]
$^2$
School of Mathematical Sciences, Key Laboratory of Intelligent Computing and Applications (Ministry of Education), Tongji University, Shanghai 200092, China, \email{2333706@tongji.edu.cn, 22210@tongji.edu.cn} \\
\end{Addresses}

\begin{abstract}
We propose a numerical approach for solving conjugate heat transfer problems using the finite volume method. This approach combines a semi-implicit scheme for fluid flow, governed by the incompressible Navier-Stokes equations, with an optimization-based approach for heat transfer across the fluid-solid interface. In the semi-implicit method, 
the convective term in the momentum equation is treated explicitly, ensuring computational efficiency, while maintaining  stability when a CFL condition involving fluid velocity is satisfied. Heat exchange between the fluid and solid domains is formulated as a constrained optimization problem, which is efficiently solved using a sequential quadratic programming method. Numerical results are presented to demonstrate the effectiveness and performance of the proposed approach.

\end{abstract}


\medskip

\section{Introduction}

Conjugate heat transfer (CHT) problems, which  involve the thermal interaction of fluid and solid domains, are prevalent in a wide range of science and engineering applications. These include
cooling system analysis \cite{srinivasan2020conjugate} and turbocharger modeling \cite{burke2015lumped} in automotive engineering, transonic nozzle flow simulation in aerospace engineering \cite{zhu2025modeling}, electronic device cooling \cite{sarma2018entropy,deb2025augmented} etc. 
Developing efficient and robust numerical methods for solving CHT problems is essential for the design, optimization, performance evaluation, and safety analysis of various industrial systems. 

Two predominant approaches are commonly employed in the numerical treatment of CHT problems:  the monolithic and partitioned methods. The monolithic approach \cite{chen2000note,pan2017fully,pan2018efficient} treats the fluid and solid governing equations as a single coupled system, offering superior convergence properties and stability. 
However, its lack of flexibility is a major drawback, as it does not allow the use of existing solvers independently for the fluid and solid sub-problems. 
In contrast, the partitioned approach 
\cite{henshaw2009composite,kazemi2014accuracy,verstraete2016stability,errera2016comparative}
separates the treatment of the fluid and solid sub-problems by solving each domain independently and iterating between them to enforce the fluid-solid interface conditions. This strategy benefits from the use of established solvers for both sub-problems, offering greater flexibility and ease of implementation.

Partitioned methods differ in their approach to enforcing the interface conditions. 
One of the most commonly used iterative techniques is the Dirichlet-to-Neumann (DtN) approach \cite{giles1997stability, kazemi2013high, verstraete2016stability, monge2018convergence}. In this method, a Dirichlet boundary condition is imposed at the fluid-solid interface for the fluid sub-problem to ensure temperature continuity, while the solid sub-problem is governed by a Neumann boundary condition at the interface to maintain heat flux continuity. Alternatively, the roles of these boundary conditions can be reversed, depending on the thermal effusivities of the fluid and solid, to ensure stability \cite{kazemi2014analysis}.
 The efficiency of partitioned schemes can be improved by exchanging boundary conditions only after multiple time steps on the fluid and solid domains, to guarantee a fast convergence of the coupled problem \cite{verstraete2016stability}. Additionally, variants such as the  Dirichlet-Robin iterative method \cite{errera2016comparative,errera2019single} have been proposed to further enhance stability and convergence, especially in fluid-solid coupling involving high Biot numbers. Furthermore, adaptive coupling coefficients have been introduced to optimize CHT for improved stability \cite{errera2025advanced}. 

In partitioned approaches for CHT problems, the fluid flow needs to be solved multiple times per time step with varying fluid-solid interface conditions due to the iterative thermal exchange between fluid and solid domains. Thus, the efficiency of the numerical solver for fluid flow, governed by the Navier-Stokes equations, is crucial.
Pressure-based methods, such as the SIMPLE-family algorithms \cite{patankar1972,darwish2016finite}, are widely used for solving incompressible fluid flows. These methods rely on a prediction-correction procedure and provide a reliable and efficient approach for CHT simulations \cite{chen2000note, oztop2008conjugate, zhu2025modeling}.
Other efficient solvers include the coupled method proposed in \cite{chen2010coupled}, which enhances computational performance by solving for pressure and velocity simultaneously, in contrast to the typical sequential approach where the momentum and pressure correction equations are solved separately. Additionally,  semi-implicit methods \cite{liles1978semi,zhang2019conservative} offer  an efficient approach to solving the Euler equations in the context of two-phase compressible flows. Furthermore, fully decoupled monolithic projection methods \cite{pan2017fully,pan2018efficient} linearize the nonlinear convective term and rely on approximate factorization techniques, ensuring  stability and computational efficiency.

The objective of this work is to propose an efficient partitioned numerical framework for solving CHT problems. To this end, we integrate a semi-implicit method for fluid flow with an optimization-based approach for  heat transfer across the fluid-solid interface. The semi-implicit method  solves the Navier-Stokes equations governing the fluid flow, treating the convective term in the momentum equation explicitly to ensure computational efficiency, while maintaining stability under a Courant-Friedrichs-Lewy (CFL) condition constrained by the fluid velocity. The stability condition is rigorously established through von Neumann analysis.  
By eliminating sub-iterations within each time step, this approach significantly enhances  computational performance for time-dependent problems. 
For heat transfer, we reformulate the coupling problem as a constrained optimization problem. While optimization-based approaches have been proposed within the finite element framework for domain decomposition problems \cite{Gunzburger_Peterson_Kwon_1999,
gunzburger2000optimization}, in this study, we extend this concept to the finite volume discretization framework. 
The resulting optimization problem is  
efficiently solved  using a sequential quadratic programming (SQP) method \cite{taddei2024non}. When  applying the SQP method,
the computation of the Jacobian matrix, which involves differentiating the sub-problem residual with respect to the control variable, is computationally expensive. This issue becomes particularly pronounced  when the number of interface faces is large. To mitigate this computational cost, we introduce a reduced representation of the interface \cite{aletti2017reduced,discacciati2023localized} to further enhance the efficiency of the optimization-based method.

 This paper is organized as follows. Section \ref{sec:governing_eqns} introduces the governing equations considered in this work. Section \ref{sec:num_methods} details the semi-implicit method for fluid flow, the optimization-based approach for heat transfer, and their integration in CHT problems. Section \ref{sec:num_results} presents numerical results to validate the proposed numerical framework and 
to illustrate its properties.

\section{Governing equations}
\label{sec:governing_eqns}

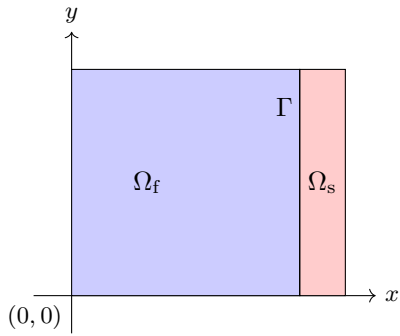
\begin{figure}[!h]
    \centering

    \begin{tikzpicture}
        \draw[->] (-0.5,0) -- (4,0) node[right] {\small $x$};
        \draw[->] (0,-0.5) -- (0,3.5) node[above] {\small $y$};
        
        \node[below left] at (0,0) {\small $(0,0)$};

        \draw [fill=blue!20] (0,0) rectangle (3,3);
        \node at (1,1.5) {${\Omega}_{\rm f}$};

        \draw [fill=red!20] (3,0) rectangle (3.6,3);
        \node at (3.3,1.5) {${\Omega}_{\rm s}$};

        \draw (3,0) -- (3,3);  
        \node at (2.8,2.5) {${\Gamma}$};  
    \end{tikzpicture}

    \caption{computational domain for conjugate heat transfer.}
    \label{fig:CHT_illustration}
\end{figure}

In this section, we present the governing equations for CHT problems 
\cite{kazemi2013high,kazemi2014accuracy,
pan2018efficient}. A schematic representation of the computational domain is shown in Figure \ref{fig:CHT_illustration}.

The fluid domain $\Omega_{\rm f}$ is governed by the 
incompressible Navier-Stokes equations for fluid flow and an advection-diffusion equation for temperature (energy equation):
\begin{align}
\label{eq:fluid_mass}
\nabla\cdot\mathbf{u}_{\rm f} &= 0, \\
\label{eq:fluid_momentum}
\rho_{\rm f}\frac{\partial \mathbf{u}_{\rm f}}{\partial t} + \nabla\cdot(\rho_{\rm f}\mathbf{u}_{\rm f}\otimes\mathbf{u}_{\rm f}) &= \nabla\cdot(\mu_{\rm f}\nabla\mathbf{u}_{\rm f})-\nabla p-
\rho_{\rm f}\beta \mathbf{g}_0(T_{\rm f}-T_{\rm {ref}}), \\
\rho_{\rm f} c_{\rm{pf}}\left(\frac{\partial T_{\rm f}}{\partial t} +  \nabla\cdot (\mathbf{u}_{\rm f}  T_{\rm f})\right) &= \nabla\cdot(k_{\rm f}\nabla T_{\rm f}) + Q_{\rm f},
\label{eq:CHT_fluid}
\end{align}
where $\mathbf{u}_{\rm f}$ is the velocity field, $p$ is the pressure, $\rho_{\rm f}$ is the fluid density, $\mu_{\rm f}$ is the dynamic viscosity, $T_{\rm f}$ is the fluid temperature, $k_{\rm f}$ is the thermal conductivity, $c_{\rm{pf}}$ is the  specific heat capacity, $Q_{\rm f}$ represents the  heat source in the fluid,   $\beta$ is the thermal expansion coefficient of the fluid, and $\mathbf{g}_0$ is the gravitational acceleration. 
The term $\rho_{\rm f}\beta \mathbf{g}_0(T_{\rm f}-T_{\rm {ref}})$ represents buoyancy effects 
under the Boussinesq approximation. If buoyancy effects are neglected, the energy equation becomes decoupled from the mass and momentum equations. 

Heat transfer in solid domain $\Omega_{\rm s}$ is governed by the heat conduction equation: 
\begin{align}
\rho_{\rm s}c_{\rm{ps}}\frac{\partial T_{\rm s}}{\partial t} = \nabla\cdot(k_{\rm s}\nabla T_{\rm s}) + Q_{\rm s},
\label{eq:CHT_solid}
\end{align}
where $ T_{\rm s} $ is the temperature field in the solid, $ k_{\rm s} $ is the thermal conductivity, $ \rho_{\rm s} $ and $ c_{\rm{ps}} $ are the density and specific heat, and $Q_{\rm s}$ is the source term. 

At the fluid-solid interface $\Gamma$, the following conjugate boundary conditions enforce continuity of  temperature and heat flux:
\begin{align}
T_{\rm s}&=T_{\rm f}, \nonumber \\
k_{\rm s}\frac{\partial T_{\rm s}}{\partial \mathbf n}&=k_{\rm f}\frac{\partial T_{\rm f}}{\partial \mathbf n},
\label{eq:CHT_coupling_conditions}
\end{align}
where $\mathbf{n}$ is the unit normal vector at the interface, pointing from the fluid domain to the solid domain. 

\section{Numerical methods}

\label{sec:num_methods}

\subsection{Semi-implicit method for the Navier-Stokes equations}
\label{sec:semi_implicit}

In this section, we present a semi-implicit finite volume method for solving the incompressible Navier-Stokes equations. This method is designed for unsteady problems, where the convective term is treated explicitly. The pressure and diffusion terms are handled implicitly to ensure that the time step is limited by the fluid flow velocity. 
By explicitly treating the nonlinear convective term, the method theoretically eliminates the need for sub-iterations to address non-linearity, unlike projection-based approaches such as the SIMPLE method. Consequently, the semi-implicit method offers computational efficiency for time-dependent problems.  
A similar semi-implicit approach has been previously considered in the context of compressible two-phase flows \cite{zhang2019conservative} using the Euler equations, where the diffusion term was not included. Here, we apply the semi-implicit method to incompressible flows and incorporate the diffusion term. 
In this study, we consider first-order time and spatial discretization. 

Here, we consider the case where the energy equation is decoupled, meaning that we solve the Navier-Stokes equations \eqref{eq:fluid_mass}-\eqref{eq:fluid_momentum} while neglecting the buoyancy term in the momentum equation \eqref{eq:fluid_momentum}. The treatment of the fully coupled system \eqref{eq:fluid_mass}-\eqref{eq:CHT_coupling_conditions}, which includes  the energy equation and the coupling with the solid, will be addressed in Section   \ref{sec:opt_cht}. 
We integrate the momentum and mass equations over a control volume $C$ using the finite volume method,  leading to the discretized Navier-Stokes equations:  
\begin{align}
\label{eq:mom_semi_implicit}
\rho_{\rm f}\frac{(\mathbf{u}_{\rm f})_{C}^{n+1}-(\mathbf{u}_{\rm f})_{C}^{n}}{\Delta t}+\frac{1}{V_{C}}\sum_{b(C)}(\rho_{\rm f}\mathbf{u}_{\rm f})_{b}^{n}\left[(\mathbf{u}_{\rm f})_{b}^{n}\cdot\mathbf{S}_{\mathrm{f}, b}\right]
&=-\nabla P_{C}^{n+1} + 
\frac{1}{V_{C}}\sum_{b(C)}(\mu_{\rm f}
\nabla\mathbf{u}_{\rm f})_{b}^{n+1}\cdot\mathbf{S}_{\mathrm f,b}, \\
\sum_{b(C)}(\mathbf{u}_{\rm f})_{b}^{n+1}\cdot\mathbf{S}_{\mathrm f,b}&=0,
\label{eq:mass_semi_implicit}
\end{align}
where $b(C)$ denotes the boundary of cell $C$, the subscript $(\cdot)_b$ indicates the value at cell boundary $b$, and the area vector $\mathbf{S}_{\mathrm{f}, b}$ is defined as $\mathbf{S}_{\mathrm{f}, b}=\mathbf{n}_{\mathrm{f}, b}A_{\mathrm{f}, b}$, i.e., the unit normal vector $\mathbf{n}_{\mathrm{f}, b}$ scaled by the area of the cell boundary $b$. 
The transported quantity $(\rho_{\rm f}\mathbf{u}_{\rm f})_{b}^{n}$ is determined using the upwind scheme, and the diffusion term $(\mu_{\rm f}
\nabla\mathbf{u}_{\rm f})_{b}^{n+1}$ can be related to the velocities of the two adjacent cells of the boundary face $b$ using standard techniques (cf. \cite[chap. 8]{darwish2016finite}). 

The discretized momentum equation \eqref{eq:mom_semi_implicit} can be rewritten in the 
following compact form
\begin{equation}
a_{C} (\mathbf{u}_{\rm f})_C^{n+1} + \sum\limits_{F\in \text{nb}({C})} a_{F} (\mathbf{u}_{\rm f})_F^{n+1} = -\nabla P_{C}^{n+1} + \mathbf{B}^n_C,
\label{eq:momentum_semi_implicit_compact}
\end{equation}
where $\text{nb}({C})$ denotes the neighboring cells of cell $C$. This compact form is exploited in  \ref{appendix:iter_solver} to derive an iterative procedure to solve the above discretized system. 

The discretized system \eqref{eq:mom_semi_implicit}-\eqref{eq:mass_semi_implicit} is further  complemented by the  Rhie-Chow interpolation \cite{rhie1983numerical}, which determines the fluid velocity at the cell face  $(\mathbf{u}_{\rm f})_{b}^{n+1}$ in the mass conservation equation \eqref{eq:mass_semi_implicit}. Specifically, the face velocity is computed as 
 \begin{equation}
 (\mathbf{u}_{\rm f})_b^{n+1}=\frac12\left[(\mathbf{u}_{\rm f})_{C}^{n+1}+(\mathbf{u}_{\rm f})_{C'}^{n+1}\right]-(\mathbf{D}_{\rm f})_b\left[(\nabla P)_b^{n+1}-\overline{(\nabla P)_b^{n+1}}\right],
 \label{eq:rhie_chow}
 \end{equation}
where $C$ and $C'$ are the two adjacent cells of face $b$; $(\mathbf{D}_{\rm f})_b=\frac12\left[(\mathbf{D}_{\rm f})_{C}+(\mathbf{D}_{\rm f})_{C'}\right]$, 
$(\mathbf{D}_{\rm f})_{C}=1/a_{C}$, $(\mathbf{D}_{\rm f})_{C'}=1/a_{C'}$; 
$\overline{(\nabla P)^{n+1}_b}$ is the average of the pressure derivatives of the two adjacent cells. In this work, we consider only orthogonal grids. For non-orthogonal grids, various correction techniques \cite{zhang2019conservative}\cite[chap. 8]{darwish2016finite}) can be employed, we omit the details for brevity.

The global system based on equations \eqref{eq:mass_semi_implicit}-\eqref{eq:rhie_chow} can be formulated as a linear system, with $(\mathbf{u}_{\rm f})_{C}^{n+1}$, $P_{C}^{n+1}$, $ (\mathbf{u}_{\rm f})_b^{n+1}$ as the unknowns. 
We present the following theorem, which establishes  that  the semi-implicit method remains stable provided that the Courant-Friedrichs-Lewy (CFL) condition, determined   by the flow velocity, is satisfied. 
\begin{theorem1}
    The semi-implicit method, with Rhie-Chow interpolation for computing the face velocity, is stable if the CFL condition is satisfied: 
\begin{equation}
    0\leqslant C_m \leqslant 1,
    \label{eq:stab_cfl}
\end{equation}
where $C_m = C_{m1} + C_{m2} = \frac{u_0\Delta t}{\Delta x}+\frac{v_0\Delta t}{\Delta y}$, $u_0$ and $v_0$ represent the characteristic flow velocities in the $x-$ and $y-$directions, respectively, $\Delta x$ and $\Delta y$ denote the spatial discretization sizes in the $x-$ and $y-$directions, and $\Delta t$ is the time step.  
\label{thm:stability}
\end{theorem1}
The proof of Theorem \ref{thm:stability} is provided in \ref{appendix:stability}. 
Since the global system is linear, standard linear solvers can be applied directly to solve it without subiterations. 
However, 
in \ref{appendix:iter_solver} we introduce a preliminary iterative method solver for the global problem, employing a prediction-correction approach similar to the SIMPLE method. Unlike the SIMPLE  method, the matrices to be solved in the prediction step and for the pressure correction equation remain unchanged across sub-iterations and time steps. This characteristic can be leveraged to further improve the efficiency of the overall algorithm. 
 The resolution of the global system \eqref{eq:mass_semi_implicit}-\eqref{eq:rhie_chow}  using more efficient techniques such as matrix factorization techniques similar to those in  \cite{pan2017fully,diaz2023non} is subject of future research.

\subsection{Optimization-based method for heat transfer problems}
\label{sec:heat_prob_opt}

We introduce an optimization-based method to solve the heat transfer problem defined by equations \eqref{eq:CHT_fluid}--\eqref{eq:CHT_solid}, subject to the coupling conditions given in \eqref{eq:CHT_coupling_conditions}. In this section, we assume that the fluid velocity, $\mathbf{u}_{\rm f}$, is known. The case where the fluid velocity is also treated as an unknown variable will be addressed in the subsequent section.

The integration of the temperature equations for the fluid and the solid over control volumes $C_{\rm f}$ and $C_{\rm s}$, respectively, yields the following equations:
\begin{align}
\label{eq:heat_disc_fluid}
\rho_{\rm f} c_{\rm{pf}} \left[\frac{T_{\rm f}^{n+1} - T_{\rm f}^{n}}{\Delta t} + \frac{1}{V_{C_{\rm f}}} \sum_{b} \left( \mathbf{u}_{\mathrm f,b} \cdot \mathbf{S}_{\mathrm f,b} \right) T_{\mathrm f,b}^{n+1} \right] &= \frac{1}{V_{C_{\rm f}}} \sum_{b} \left( k_{\mathrm f,b} \nabla T_{\mathrm f,b}^{n+1} \cdot \mathbf{S}_{\mathrm f,b} \right) + Q_{\rm f}, \\
\rho_{\rm s} c_{\rm{ps}} \frac{T_{\rm s}^{n+1} - T_{\rm s}^{n}}{\Delta t} &= \frac{1}{V_{C_{\rm s}}} \sum_{b} \left( k_{\mathrm s,b} \nabla T_{\mathrm s,b}^{n+1} \cdot \mathbf{S}_{\mathrm s,b} \right) + Q_{\rm s}.
\end{align}
Here, the area vector $\mathbf{S}_{\mathrm f,b}$ is defined as in Section \ref{sec:semi_implicit}, and $\mathbf{S}_{\mathrm s,b}$ is defined in a similar manner.  
We distinguish between two types of cell boundaries: $(1)$ internal boundaries or domain boundaries that do not involve the fluid-solid interface, and $(2)$ boundaries corresponding to the fluid-solid interface. For the first type of cell boundary, standard numerical strategies can be applied to handle the associated terms. Specifically, the convective term in the fluid equation at the cell boundary $T_{\mathrm f,b}^{n+1}$ can be determined using an upwind scheme or by applying appropriate boundary conditions. The diffusion terms at the cell boundary, $\nabla T_{\mathrm f,b}^{n+1}$ and $\nabla T_{\mathrm s,b}^{n+1}$, can be related to the temperatures of the two adjacent cells at the boundary face $b$ using standard techniques (cf. \cite[chap. 8]{darwish2016finite}), or by applying the appropriate treatment of boundary conditions.

\begin{figure}[!h]

\centering
\begin{tikzpicture}[scale=0.8]

\draw[fill=blue!20] (0,0) rectangle (3,3);
\node at (1.2, 1.5) {$C_{\rm f}$};

\filldraw[black] (1.5, 1.5) circle (2pt);  

\draw[fill=red!20] (3,0) rectangle (6,3);
\node at (4.2, 1.5) {$C_{\rm s}$};

\filldraw[black] (4.5, 1.5) circle (2pt);  

\draw[thick] (3,-1) -- (3,4);  
\node at (3.3, 1.5) {$b$};

\node at (3, 1.5) {\textbf{+}};  

\draw[->, thick, blue] (3, 2.5) -- (4, 2.5);  
\node at (4.4, 2.5) {$\mathbf{n}_{\mathrm{f},b}$};  

\draw[->, thick, red] (3, 0.5) -- (2, 0.5);  
\node at (1.6, 0.5) {$\mathbf{n}_{\mathrm{s},b}$};  

\end{tikzpicture}
\caption{illustration of fluid-solid interface with two adjacent cells. }
\label{fig:Tdiff_interface}
\end{figure}
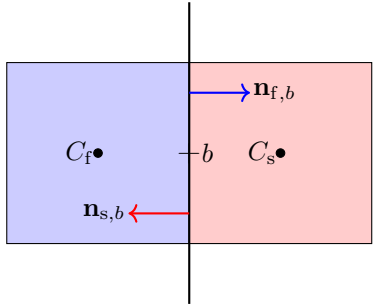

For a cell boundary of type $(2)$, i.e., a boundary that lies on the fluid-solid interface (cf. Figure \ref{fig:Tdiff_interface}, where the cell boundary $b$ located at the interface), the contribution of the cell boundary to the convective term in the fluid equation is zero, as we assume a no-slip wall condition at the interface.
Regarding the diffusion term, the coupling condition \eqref{eq:CHT_coupling_conditions}$_2$ and the assumption of a conforming mesh ensure that the heat flux remains continuous across the interface, leading to the equality $k_{\mathrm f,b} \nabla T_{\mathrm f,b}^{n+1} \cdot \mathbf{S}_{\mathrm f,b} = -k_{\mathrm s,b} \nabla T_{\mathrm s,b}^{n+1} \cdot \mathbf{S}_{\mathrm s,b}$. We denote this common heat flux across the interface face $b$ as $g_b^{n+1}$, such that:  
\begin{equation}
g_b^{n+1}= k_{\mathrm f,b} \nabla T_{\mathrm f,b}^{n+1} \cdot \mathbf{S}_{\mathrm f,b} = -k_{\mathrm s,b} \nabla T_{\mathrm s,b}^{n+1} \cdot \mathbf{S}_{\mathrm s,b}.
\label{eq:g_def}
\end{equation}
This formulation leads to the following algebraic equations in compact form for the heat transfer equations in the fluid and solid subdomains:
\begin{align}
\label{eq:opt_based_fluid}
\mathbf{R}_{\rm f}(\mathbf{T}_{\rm f}^{n+1}) - \mathbf{E}_{\rm f} \mathbf{g}^{n+1} = 0, \\
\mathbf{R}_{\rm s}(\mathbf{T}_{\rm s}^{n+1}) + \mathbf{E}_{\rm s} \mathbf{g}^{n+1} = 0,
\label{eq:opt_based_solid}
\end{align}
where $\mathbf{T}_{\rm f}^{n+1} \in \mathbb{R}^{N_{\rm f}}$ and $\mathbf{T}_{\rm s}^{n+1} \in \mathbb{R}^{N_{\rm s}}$ represent the temperature vectors at the cell center for all cells in the fluid and solid domains, respectively. Here, $N_{\rm f}$ and $N_{\rm s}$ denote the number of cells in the fluid and solid domains. The residuals $\mathbf{R}_{\rm f} \in \mathbb{R}^{N_{\rm f}}$ and $\mathbf{R}_{\rm s} \in \mathbb{R}^{N_{\rm s}}$ may be nonlinear, for example, when temperature-dependent diffusion coefficients $k_{\rm f}$ and $k_{\rm s}$ are considered. 
The vector $\mathbf{g}\in \mathbb{R}^{N_{\Gamma}}$ represents the heat flux across all faces of the interface, where $N_{\Gamma}$ is the total number of interface faces. 
In accordance with the discretized equation \eqref{eq:heat_disc_fluid}, the matrix $\mathbf{E}_{\rm f} \in \mathbb{R}^{N_{\rm f} \times N_{\Gamma}}$  in the coupling term is constructed as follows: $(\mathbf{E}_{\rm f})_{ij} = 1$ if the $j$-th face belongs to the $i$-th cell, and $(\mathbf{E}_{\rm f})_{ij} = 0$ otherwise. The construction for $\mathbf{E}_{\rm s} \in \mathbb{R}^{N_{\rm s} \times N_{\Gamma}}$ follows similarly.

Note that given the flux at the fluid-structure interface $\mathbf{g}$, the fluid and solid sub-problems can be decoupled and fully determined independently. As described in \cite{Gunzburger_Peterson_Kwon_1999,
gunzburger2000optimization,taddei2024non}, we refer to $\mathbf{g}$ as the control. 
However, an arbitrary choice of the control $\mathbf{g}$ does not guarantee that the local solutions to \eqref{eq:opt_based_fluid}-\eqref{eq:opt_based_solid}
satisfy the transmission condition \eqref{eq:CHT_coupling_conditions}$_1$ and hence will not be
valid solutions to the coupled problem \eqref{eq:CHT_fluid} and \eqref{eq:CHT_solid}. 
The \emph{optimal control} $\mathbf{g}$
that guarantees temperature continuity at the interface $\Gamma$ can be obtained by solving the following optimization problem
\begin{equation}
\min_{\substack{
\mathbf{T}_{\rm f}^{n+1} ; \\
\mathbf{T}_{\rm s}^{n+1} ; \\
\mathbf{g}^{n+1} } }\; 
F_{\delta}(\mathbf{T}_{\rm f}^{n+1},\mathbf{T}_{\rm s}^{n+1},\mathbf{g}^{n+1}):= 
 \frac{1}{2}\left\Vert \mathbf{T}_{\rm f,\Gamma}^{n+1}-\mathbf{T}_{\rm s,\Gamma}^{n+1}\right\Vert_2^2 + \frac{\delta}{2}\Vert \mathbf{g}^{n+1}\Vert_2^2,
 \label{eq:opt_statement}
\end{equation}
with equations \eqref{eq:opt_based_fluid} and \eqref{eq:opt_based_solid} as constraints. Here, $\mathbf{T}_{\rm f,\Gamma}^{n+1}$ and $\mathbf{T}_{\rm s,\Gamma}^{n+1}$  represent the temperature at all faces of the interface $\Gamma$ for the fluid and solid domains, respectively. 
The second term in the objective function of \eqref{eq:opt_statement} is a regularizer, designed to penalize controls of excessive size; the positive constant $\delta$ is selected to control the relative importance of the two terms in the objective function.
The optimization-based methodology has been effectively applied to domain-decomposition problems  involving elliptic equations \cite{Gunzburger_Lee_2000_elliptic} and Navier-Stokes equations \cite{gunzburger2000optimization, taddei2024non} within the framework of finite-element methods. In this work, we extend this methodology   to the context of finite-volume methods for addressing multiphysics problems. 

In the optimization formula, the temperature at the interface is required in the objective function \eqref{eq:opt_statement}. However, in the context of finite volume methods, the resolution of the local subproblems \eqref{eq:opt_based_fluid}-\eqref{eq:opt_based_solid} gives the temperature at the cell center rather than directly at the cell face. We now describe how the temperature at the  interface can be determined using the heat flux $\mathbf{g}$. The directional gradient of the temperature along the normal at the interface boundary $b$ (cf. Figure \ref{fig:Tdiff_interface}) can be determined as
\begin{equation*}
\nabla T_{\mathrm f,b}^{n+1}\cdot \mathbf{n}_{\mathrm f,b} = \frac{ g_b^{n+1} }{k_{\mathrm f,b} {A}_{\mathrm f,b}},
\end{equation*}
for the fluid, and similarly for the solid
\begin{equation*}
\nabla T_{\mathrm s,b}^{n+1}\cdot \mathbf{n}_{\mathrm s,b} = \frac{ -g_b^{n+1} }{k_{\mathrm s,b} {A}_{\mathrm s,b}},
\end{equation*}
where these expressions follows directly from the definition of $g_b^{n+1}$ in equation \eqref{eq:g_def}. The temperature at the interface can now be determined using these gradients and the temperature at the cell center as
\begin{equation*}
T_{\mathrm f,b}^{n+1} = \nabla T_{\mathrm f,b}^{n+1}\cdot \mathbf{n}_{\mathrm f,b} d_{b, C_{\rm f}} + T_{C_{\rm f}}^{n+1},
\end{equation*} 
where $d_{b, C_{\rm f}}$ is the distance between the cell center $C_{\rm f}$ and the face center $b$. 
This relation can be established for all faces of the interface, and finally we   
 obtain the temperature at the interface in vector form 
 \begin{equation}
 \mathbf{T}^{n+1}_{\rm f,\Gamma} = \frac{\mathbf{g}^{n+1}}{\mathbf{k}_{\mathrm f,\Gamma}\mathbf{A}_{\mathrm f,\Gamma}/\mathbf{d}_{\mathrm f,\Gamma}} + \mathbf{T}_{\rm f}^{n+1}\vert_{\Gamma},
 \label{eq:temp_interface_center_f}
 \end{equation}
 where $\mathbf{k}_{\mathrm f,\Gamma}$, $\mathbf{A}_{\mathrm f,\Gamma}$, $\mathbf{d}_{\mathrm f,\Gamma}$ represent the vectorized thermal conductivity, face area, and distance between the face center and the corresponding cell center, for all faces of the interface;  $\mathbf{T}_{\rm f}^{n+1}\vert_{\Gamma}$ denotes the restriction of $\mathbf{T}_{\rm f}^{n+1}$ to cells that share a face at the interface. A similar expression can be obtained for the solid:
  \begin{equation}
 \mathbf{T}^{n+1}_{\rm s,\Gamma} = \frac{\mathbf{g}^{n+1}}{\mathbf{k}_{\mathrm s,\Gamma}\mathbf{A}_{\mathrm s,\Gamma}/\mathbf{d}_{\mathrm s,\Gamma}} + \mathbf{T}_{\rm s}^{n+1}\vert_{\Gamma}.
 \label{eq:temp_interface_center_s}
 \end{equation}

 \subsubsection{Sequential quadratic programming method for solving the optimization problem}
\label{sec:sqp} 
 
 We resort to the sequential quadratic programming (SQP) method \cite{taddei2024non} to solve the constrained optimization problem, which involves the linearization of the nonlinear constraints. Given the relations \eqref{eq:temp_interface_center_f}-\eqref{eq:temp_interface_center_s}, which link the temperature at cell face and the temperature at cell center, we solve  the following optimization problem at each SQP iteration: 
\begin{align}
\min_{\substack{
\mathbf{T}_{\rm f},
\mathbf{T}_{\rm s},
\mathbf{g} } }\; 
F_{\delta}(\mathbf{T}_{\rm f},\mathbf{T}_{\rm s},\mathbf{g}):= 
 \frac{1}{2}\left\Vert
 (\mathbf{T}_{\rm f}
 -
 \mathbf{T}_{\rm s})\vert_{\Gamma}
 +
\mathbf{g}\left( \frac{1}{\mathbf{k}_{\rm {f,\Gamma}}\mathbf{A}_{\rm {f,\Gamma}}/\mathbf{d}_{\rm {f,\Gamma}}}
  +
  \frac{1}{\mathbf{k}_{\rm {s,\Gamma}}\mathbf{A}_{\rm {s,\Gamma}}/\mathbf{d}_{\rm {s,\Gamma}}} \right)
  \right\Vert_2^2 + \frac{\delta}{2}\Vert \mathbf{g}\Vert_2^2,
   \label{eq:opt_statement_sqp}
\end{align} 
subject to the linearized constraints
 \begin{equation}
 \left\{
 \begin{array}{l}
 \mathbf{R}_{\rm f}^{n+1,k} + \mathbf{J}_{\rm f}^{{n+1,k}}(\mathbf{T}_{\rm f} - \mathbf{T}_{\rm f}^{n+1,k}) - \mathbf{E}_{\rm f}\mathbf{g} = 0, \\
\mathbf{R}_{\rm s}^{n+1,k} + \mathbf{J}_{\rm s}^{{n+1,k}}(\mathbf{T}_{\rm s} - \mathbf{T}_{\rm s}^{n+1,k})  + \mathbf{E}_{\rm s}\mathbf{g} = 0,
 \end{array}
\right.
\label{eq:opt_based_solid_lin}
\end{equation}
where $k$ denotes the SQP iteration number,
 $\mathbf{R}_{\rm f}^{n+1,k}$ and $\mathbf{R}_{\rm s}^{n+1,k}$ are the residuals evaluated at
$\mathbf{T}_{\rm s}^{n+1,k}$, $\mathbf{T}_{\rm s}^{n+1,k}$, respectively; the Jacobian matrices
 $\mathbf{J}_{\rm f}^{{n+1,k}}$ and $\mathbf{J}_{\rm s}^{{n+1,k}}$ represent the derivatives of the residuals $\mathbf{R}_{\rm f}^{n+1,k}$ and $\mathbf{R}_{\rm s}^{n+1,k}$ with respect to the temperature evaluated at iteration $k$, $\mathbf{T}_{\rm s}^{n+1,k}$, $\mathbf{T}_{\rm s}^{n+1,k}$. The linearized constraints lead to a relationship between $\mathbf{T}_{\rm  f}$, $\mathbf{T}_{\rm  s}$ and the control $\mathbf{g}$. By substituting these relations in the objective function \eqref{eq:opt_statement_sqp}, we obtain a least-squares problem for the control $\mathbf{g}$, which can be solved efficiently using the Gauss-Newton method \cite[chap. 10]{nocedal2006numerical}. The SQP method converges as long as the difference in $\mathbf{g}$ between two successive SQP iterations is below a predefined threshold,  which allows to proceed to the next time step. 
 
 \begin{remark}
One of the most commonly used approaches for solving CHT problems is the Dirichlet-to-Neumann (DtN) method \cite{verstraete2016stability}. In this method, the temperature is computed on the solid side of the interface and used as a Dirichlet boundary condition for the fluid domain;  subsequently, the heat flux is calculated on the fluid side and applied as a Neumann boundary condition on the solid side; or vice versa. This procedure is inherently sequential and relies on alternating boundary conditions between the two subdomains.
In contrast, the optimization-based (OB) method offers a parallel approach by introducing a control variable, $\mathbf{g}$, which represents the heat flow across the interface. In this formulation, 
both subdomains are treated with  Neumann boundary conditions.  The coupling condition for the continuity of heat flux is inherently satisfied, while the continuity of temperature is ensured by solving the optimization problem \eqref{eq:opt_statement}.
In Section \ref{sec:num_results}, 
we compare the performance of the DtN method and the OB method, highlighting the efficiency of the OB method. 
 
 \end{remark}

\subsubsection{Reduced optimization-based method}
\label{sec:reduced_OB}

In applying the SQP method to solve the optimization problem, each SQP iteration requires the computation of
$\left(\mathbf{J}_{\rm f}^{{n+1,k}}\right)^{-1}\mathbf{E}_{\rm f}$ and $\left(\mathbf{J}_{\rm s}^{{n+1,k}}\right)^{-1}\mathbf{E}_{\rm s}$, which  can be computationally expensive, 
especially when the number of faces at the interface is large. To improve the efficiency of the SQP method, we adopt a reduced representation of the control as follows: 
\begin{equation}
\mathbf{g} = \boldsymbol{\Phi\beta},
\label{eq:reduced_g}
\end{equation}
where $\boldsymbol{\Phi}=(\phi_1,\phi_2,\cdots,\phi_{N^r_{\Gamma}})^{\top} \in \mathbb{R}^{N_{\Gamma}\times N^r_{\Gamma}}$ represents the basis functions of the reduced space, $N^r_{\Gamma}$ is the dimensionality of the reduced space, the vector $\boldsymbol{\beta}\in \mathbb{R}^{N^r_{\Gamma}} $ is the generalized coordinates of $\mathbf{g}$ within the reduced space, which are to be determined. 

We follow the approach presented in 
\cite{aletti2017reduced,discacciati2023localized} to define the basis functions $\phi_j$ as eigenfunctions of the Laplace-Beltrami operator on the interface $\Gamma$. Specifically, these eigenfunctions satisfy 
\begin{equation}
-\Delta \phi_j = \mu_j \phi_j,
\end{equation}
where $\mu_j$ denotes the eigenvalue corresponding to the eigenfunction $\phi_j$. 
The advantages of using eigenfunctions of Laplace-Beltrami operator are detailed in \cite{aletti2017reduced}. For the problems considered in this work, where the interface is one-dimensional, the eigenfunctions $\phi_j$ reduce to a set of trigonometric functions \cite{discacciati2023localized}, given by 
$$\phi_j\in\left\{1,\:\ldots\:,\:\cos\!\left(\frac{(N^{r}-1)\pi}L\xi\right),\:\sin\!\left(\frac{\pi}L\xi\right)\:\ldots\:,\:\sin\!\left(\frac{(N^{r}-1)\pi}L\xi\right)\right\},$$
where $\xi$ is a curvilinear coordinate defined on the interface and $L$ the total length of the interface, $N^r$ is the number of retained trigonometric functions.  Note that $N^r_{\Gamma} = 2N^r-1$. 

By substituting the reduced representation of the control $\mathbf{g}$ (as given in equation  \eqref{eq:reduced_g}) into the objective function \eqref{eq:opt_statement_sqp} and the linearized constraints \eqref{eq:opt_based_solid_lin}, we reformulate the optimization problem  to solve for $(\mathbf{T}_{\rm f},\mathbf{T}_{\rm s},\boldsymbol{\beta})$. To address this optimization problem, as outlined in the previous section, we first derive the relationship between $\mathbf{T}_{\rm f},\mathbf{T}_{\rm s}$ and $\boldsymbol{\beta}$ from the linearized constraints. These relations are then substituted into the objective function to formulate a least-squares problem for the generalized coordinates 
 $\boldsymbol{\beta}$ at each SQP iteration. 
 We remark that to express $\mathbf{T}_{\rm f},\mathbf{T}_{\rm s}$ as function of $\boldsymbol{\beta}$, we now need to compute $\left(\mathbf{J}_{\rm f}^{{n+1,k}}\right)^{-1}(\mathbf{E}_{\rm f}\boldsymbol{\Phi})$ and $\left(\mathbf{J}_{\rm s}^{{n+1,k}}\right)^{-1}(\mathbf{E}_{\rm s}\boldsymbol{\Phi})$, which is significantly more efficient compared to the full computation when $N^r_{\Gamma} \ll N_{\Gamma} $. 
Hereafter, we refer to this approach as the 
reduced optimization-based method.

\subsection{Optimization-based method for conjugate heat transfer problems}
\label{sec:opt_cht}

In Algorithm \ref{algo:OB_CHT}, we outline the extension of the optimization-based method for solving CHT problems. The algorithm demonstrates the procedure for a single time step, with line $2$ representing the SQP iteration loop, where $k$ denotes the SQP iteration number. 
The integration of the fluid solver into the SQP loop is achieved in line $3$. This step ensures the velocity field from the fluid domain is updated iteratively to account for the interaction between fluid dynamics and heat transfer. Note that if we assume the fluid velocity is known (i.e., line $3$ is removed), Algorithm \ref{algo:OB_CHT} reduces to the SQP method for heat transfer problems alone, as described in Section \ref{sec:sqp}, where only the heat transfer problem is iteratively solved.

\begin{algorithm}[H]                      
\caption{optimization-based algorithm for CHT problems.}

\begin{algorithmic}[1]
\State Initialize $\mathbf{T}_{\rm f}^{n+1,0},\,\mathbf{g}^{n+1,0}$.
\vspace{3pt}
\For {$k=0, \ldots, K$ }
    \State Solve the fluid flow equations \eqref{eq:fluid_mass}-\eqref{eq:fluid_momentum} 
     using the temperature field $ \mathbf{T}_{\rm f}^{n+1,k}$
     in the buoyancy term. 
\Statex \hspace{1.5em}  This step yields the velocity field $\mathbf{u}_{\rm f}^{n+1,k}$. 

    \State Solve the optimization problem \eqref{eq:opt_statement_sqp}-\eqref{eq:opt_based_solid_lin} for
     the heat transfer problem 
    \Statex \hspace{1.5em}  using $\mathbf{u}_{\rm f} = \mathbf{u}_{\rm f}^{n+1,k}$ for the convective term
     in the fluid heat equation. 

    \Statex \hspace{1.5em} This step provides updated values for $\mathbf{T}_{\rm f}^{n+1,k+1},\,\mathbf{T}_{\rm s}^{n+1,k+1}$, and $\mathbf{g}^{n+1,k+1}$.
    \State Terminate if the convergence criterion 
    $
    \frac{\Vert \mathbf{g}^{n+1,k+1} - \mathbf{g}^{n+1,k} \Vert}{\Vert \mathbf{g}^{n+1,k} \Vert} < \texttt{tol}
    $
    is satisfied.
    \State Update the variables as follows:
    $
    \mathbf{T}_{\rm f}^{n+1,k} \gets \mathbf{T}_{\rm f}^{n+1,k+1}, \quad
    \mathbf{T}_{\rm s}^{n+1,k} \gets \mathbf{T}_{\rm s}^{n+1,k+1}, \quad
    \mathbf{g}^{n+1,k} \gets \mathbf{g}^{n+1,k+1}.
    $
\EndFor
\end{algorithmic}
\label{algo:OB_CHT}
\end{algorithm}

\begin{remark}
The optimization-based method is highly versatile and can be combined with any fluid solver, such as the SIMPLE method. To adapt the algorithm to a specific solver, it is sufficient to replace the step in line $3$ with the corresponding method for solving the fluid flow.
\end{remark}

\section{Numerical results}
\label{sec:num_results}

We present numerical results obtained using the semi-implicit method and the optimization-based method to illustrate overall numerical properties. The lid-driven cavity case illustrate the performance of the semi-implicit method; while the diffusion problem show the effectiveness of the OB method. Additionally, two CHT problems are simulated to further evaluate the methods. All numerical simulations are performed in \texttt{MATLAB\_R2024a} on a commodity Laptop. 

\subsection{Lid-driven cavity}
\begin{figure}[H]
\centering

\subfloat[$\text{Re}=100$]{\includegraphics[width=0.32\textwidth]{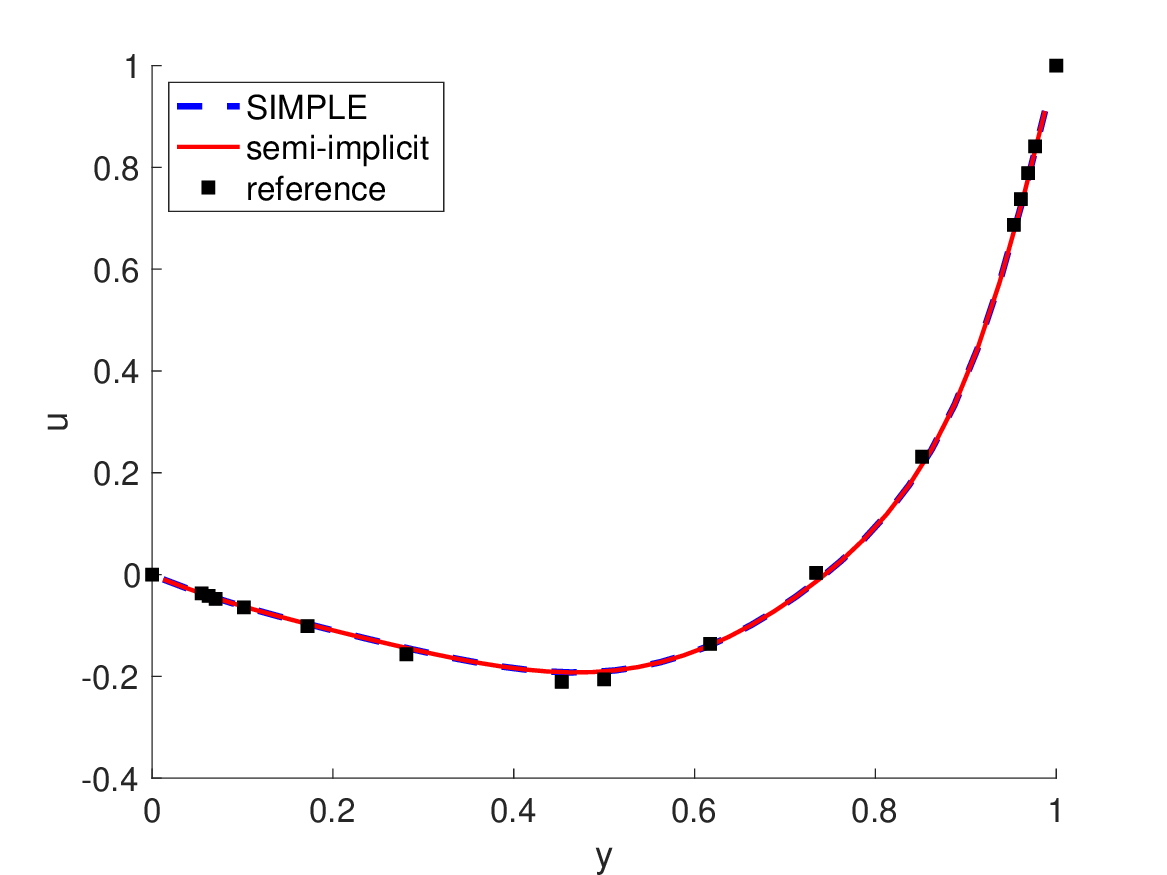}}
~~
\subfloat[$\text{Re}=400$]{\includegraphics[width=0.32\textwidth]{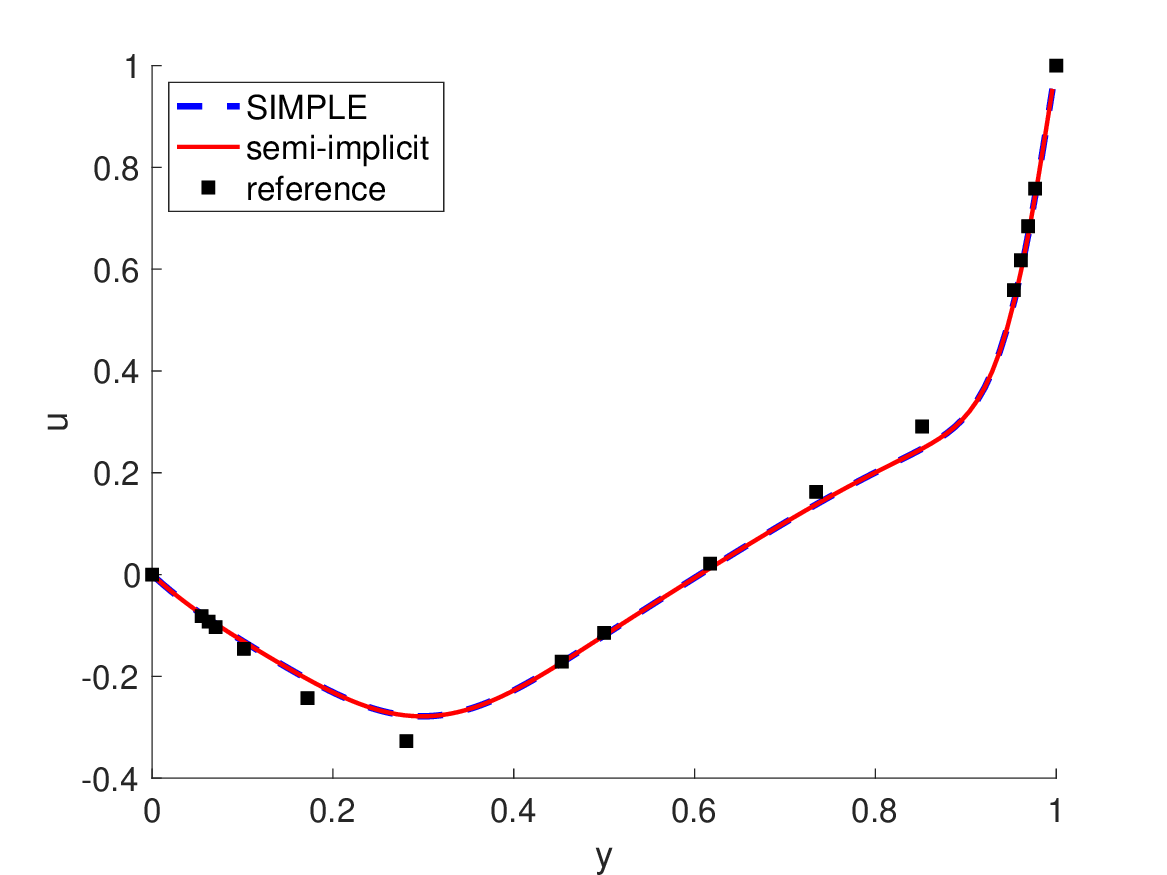}}
~~
\subfloat[$\text{Re}=1000$]{\includegraphics[width=0.32\textwidth]{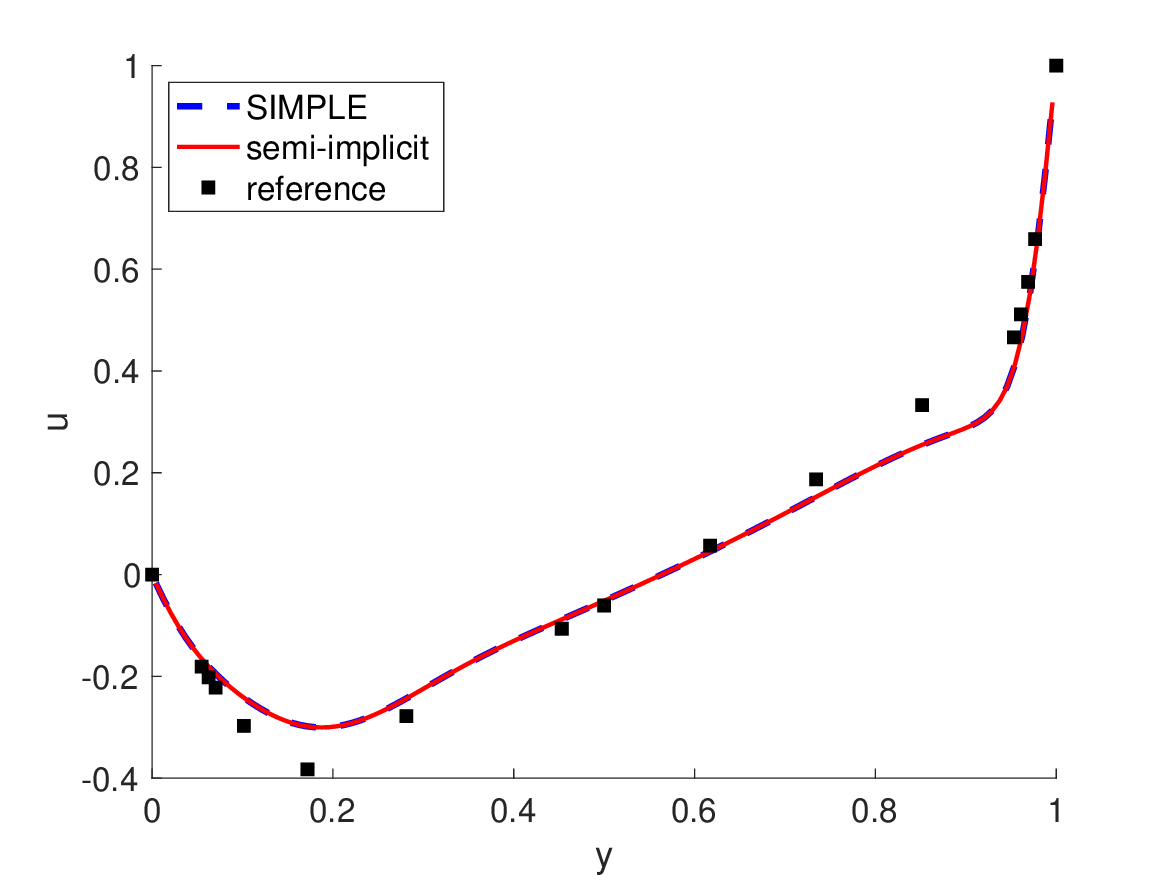}}
~~

\subfloat[$\text{Re}=100$]{\includegraphics[width=0.32\textwidth]{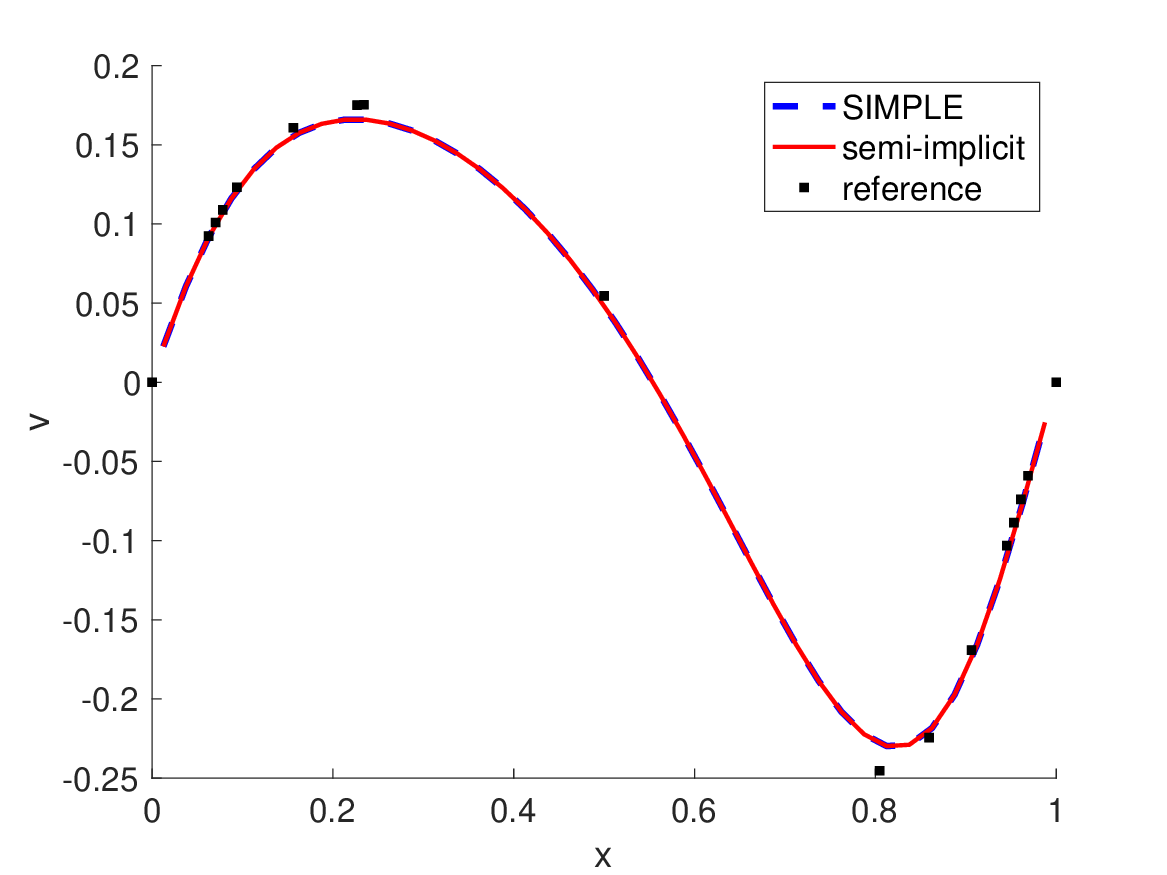}}
~~
\subfloat[$\text{Re}=400$]{\includegraphics[width=0.32\textwidth]{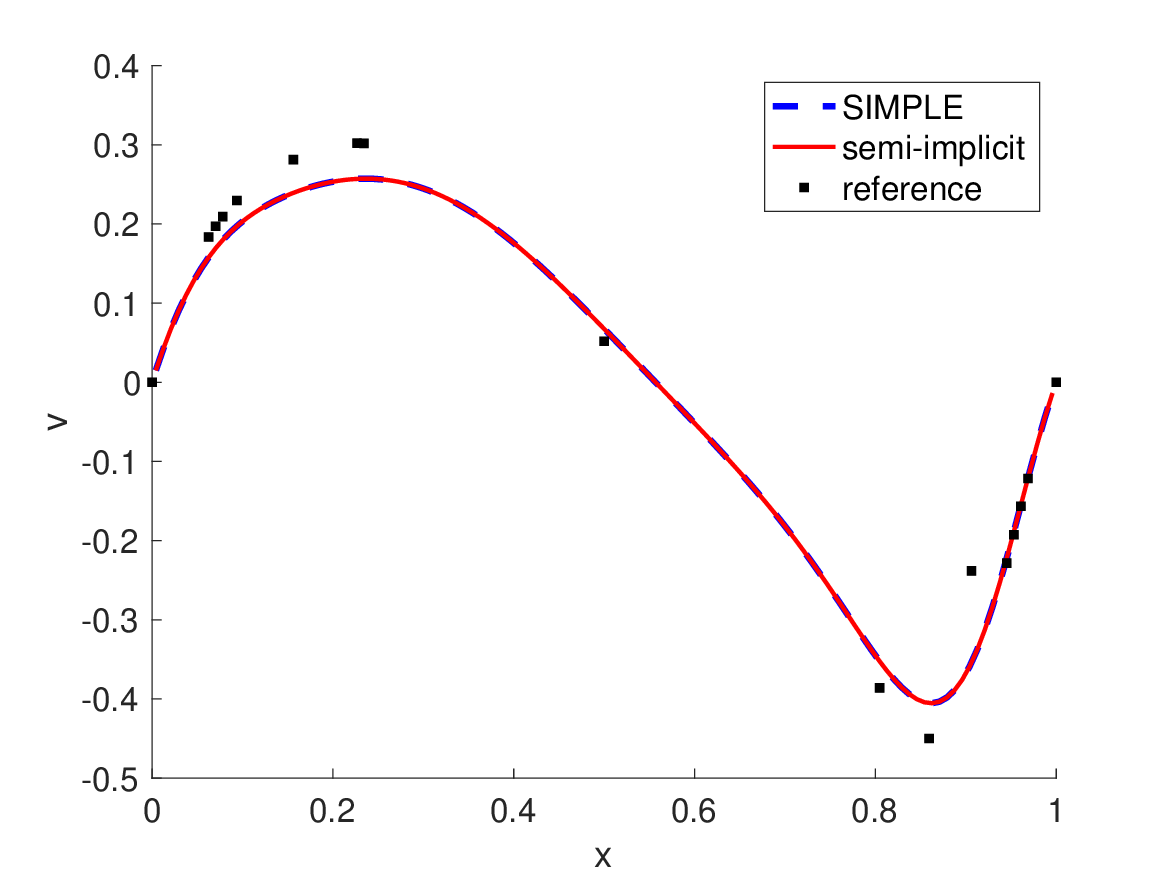}}
~~
\subfloat[$\text{Re}=1000$]{\includegraphics[width=0.32\textwidth]{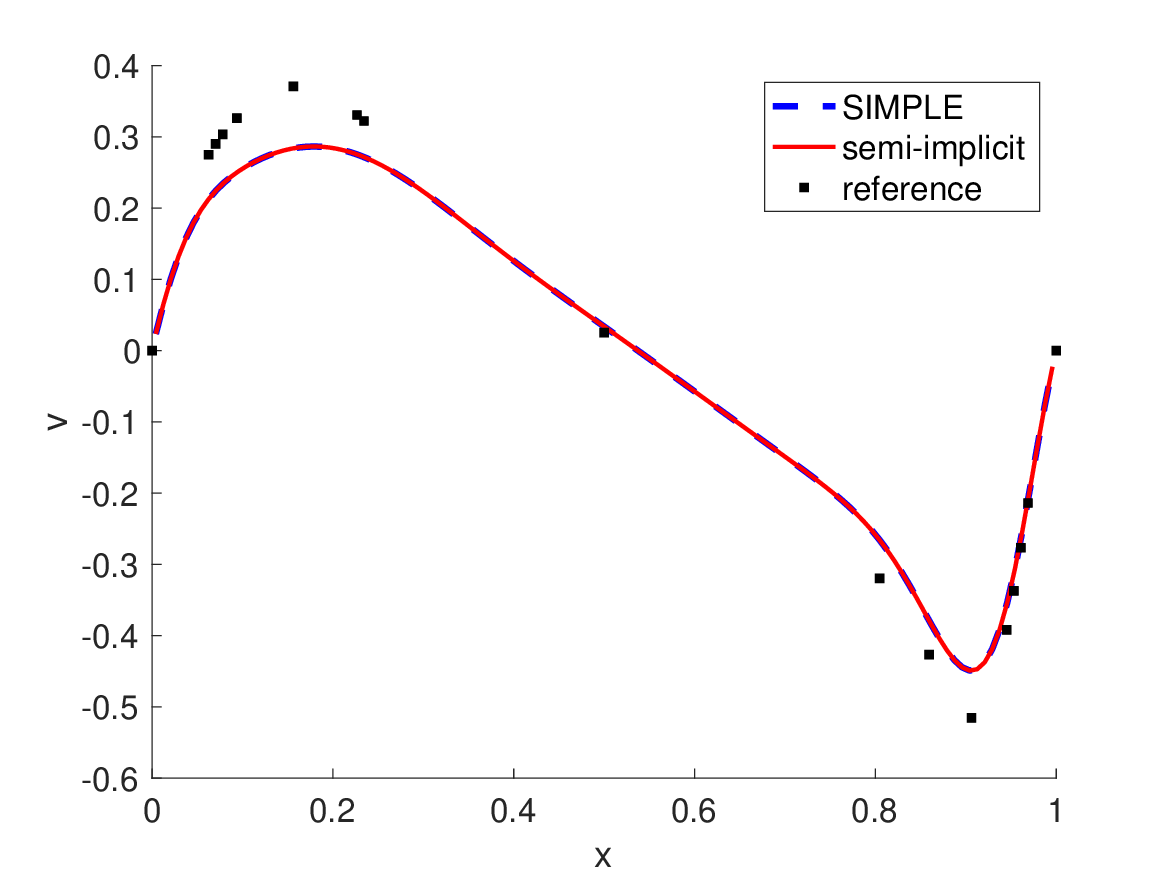}}
\caption{lid-driven cavity: (a), (b), (c) $x$-velocity  along vertical line through geometric center of cavity; (d), (e), (f) $y$-velocity along horizontal line through geometric center of cavity.}
\label{fig:res_cavity_u_slice}
\end{figure}

\begin{figure}[H]
\centering
{\includegraphics[width=8cm]{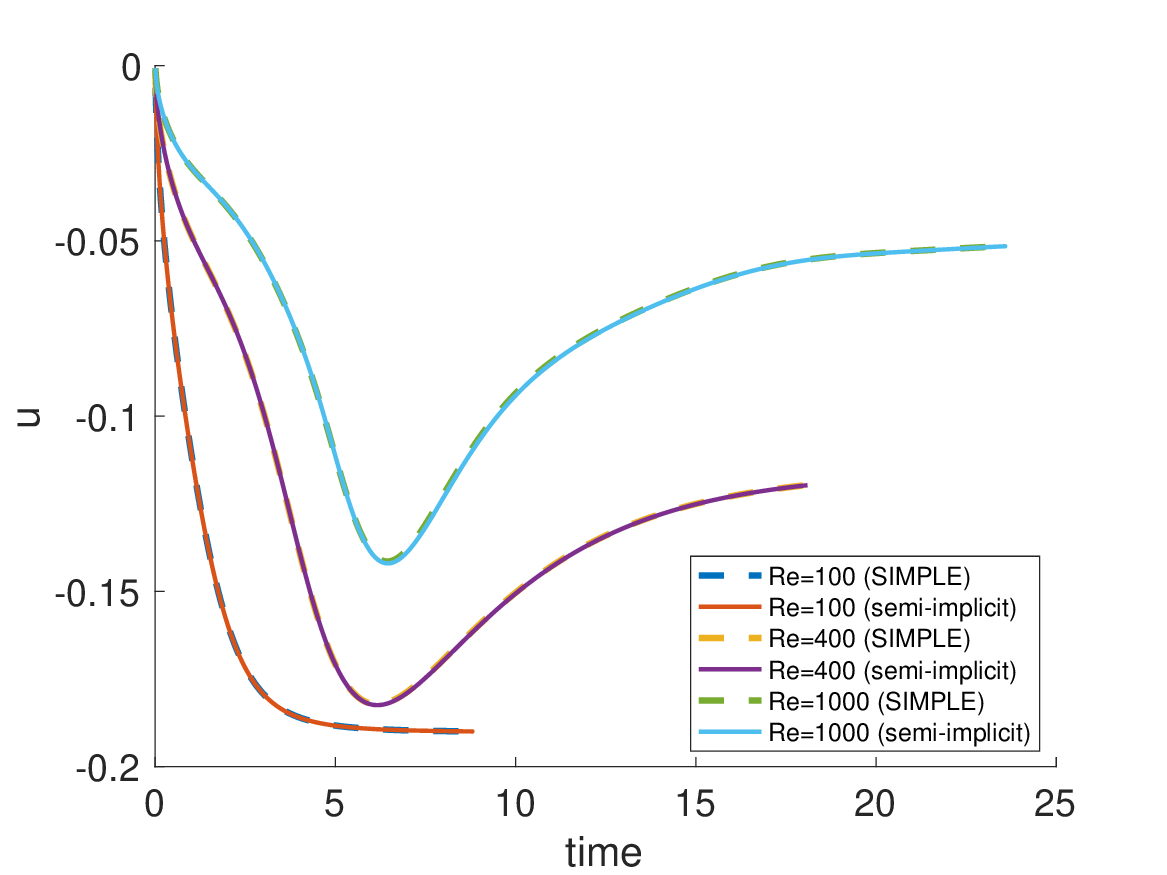}}
\caption{lid-driven cavity: evolution of the $x$-velocity at the center point.}
\label{fig:res_cavity_u_probe}
\end{figure}

We apply the semi-implicit method to simulate laminar, incompressible flow in a square cavity with a side length of $L = 1$ \cite{ghia1982high, chen2010coupled}. The left, right, and bottom boundaries are treated as no-slip walls, while the upper boundary is a moving top lid with a horizontal velocity of $u_0 = 1$. We vary the fluid viscosity $\mu$ to investigate flow behavior at different Reynolds numbers, $\text{Re} = \frac{\rho u_0 L}{\mu} = 100, 400, 1000$, with a fluid density of $\rho = 1$. Structured grids with $N \times N$ cells are used, where $N = 40$ for $\text{Re} = 100$ and $N = 120$ for $\text{Re} = 400$ and $\text{Re} = 1000$. The time step is set to $\Delta t = \text{CFL} \frac{\Delta x}{u_0}$, with a CFL number of $0.5$, where $\Delta x$ represents the characteristic mesh size. 
For this test case, no iterations are performed 
at each time step in the semi-implicit method (i.e., we set $K = 0$ for the iterative solver described in \ref{appendix:iter_solver}). 
In Figure \ref{fig:res_cavity_u}, we present the $x$-velocity field and the streamline pattern obtained  at steady state using both the semi-implicit and SIMPLE methods for $\text{Re} = 1000$. Both methods yield comparable results. The steady state is achieved when the difference between two successive velocity fields is less than $10^{-5}$: 
\begin{equation}
\text{err}=\frac{\sqrt{\sum_{i=1}^{\mathrm{Ncell}}\vert\mathbf{u}^{n+1}_i-\mathbf{u}^n_i\vert^2}}{\sqrt{\sum_{i=1}^{\mathrm{Ncell}}\vert\mathbf{u}^{n+1}_i\vert^2}}<10^{-5}.
\label{eq:steady_state}
\end{equation}
 In Figure \ref{fig:res_cavity_u_slice}, we present the velocity profile along the vertical line passing through the center of the cavity for different Reynolds numbers, compared to the reference results from \cite{ghia1982high}. The results show good agreement with the reference data. In Figure \ref{fig:res_cavity_u_probe}, we compare the time evolution of the velocity at the center point, using both the semi-implicit and SIMPLE methods. The results from both methods align well with each other.  
During the initial stages of the simulation, the SIMPLE method requires approximately $10$ iterations per time step to converge, whereas the semi-implicit method, which does not rely on an iterative procedure for this test case, demonstrates superior computational efficiency. Despite the absence of an iterative procedure, the semi-implicit is still able to capture accurately the unsteady behavoir of the flow. 

\begin{figure}[H]
\centering

\subfloat[semi-implicit method]{\includegraphics[width=0.45\textwidth]{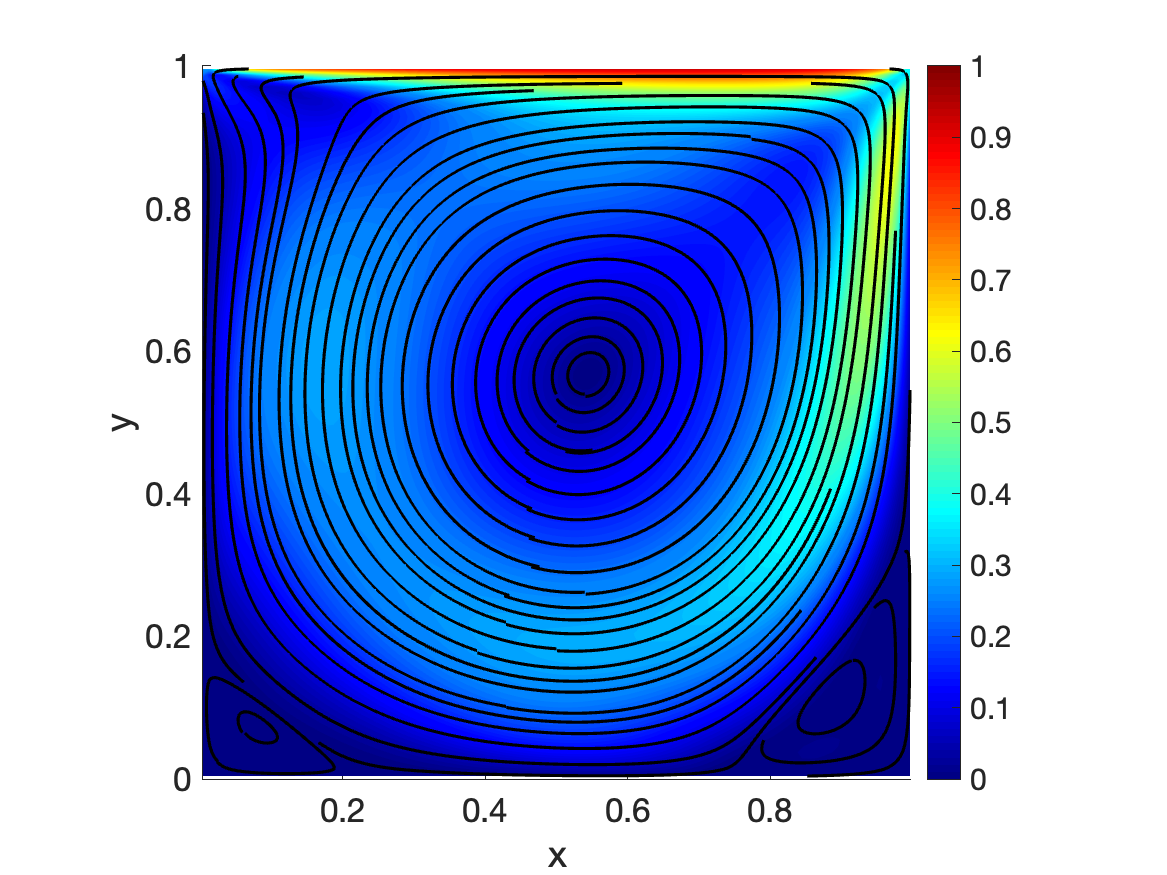}}
~~
\subfloat[SIMPLE method]{\includegraphics[width=0.45\textwidth]{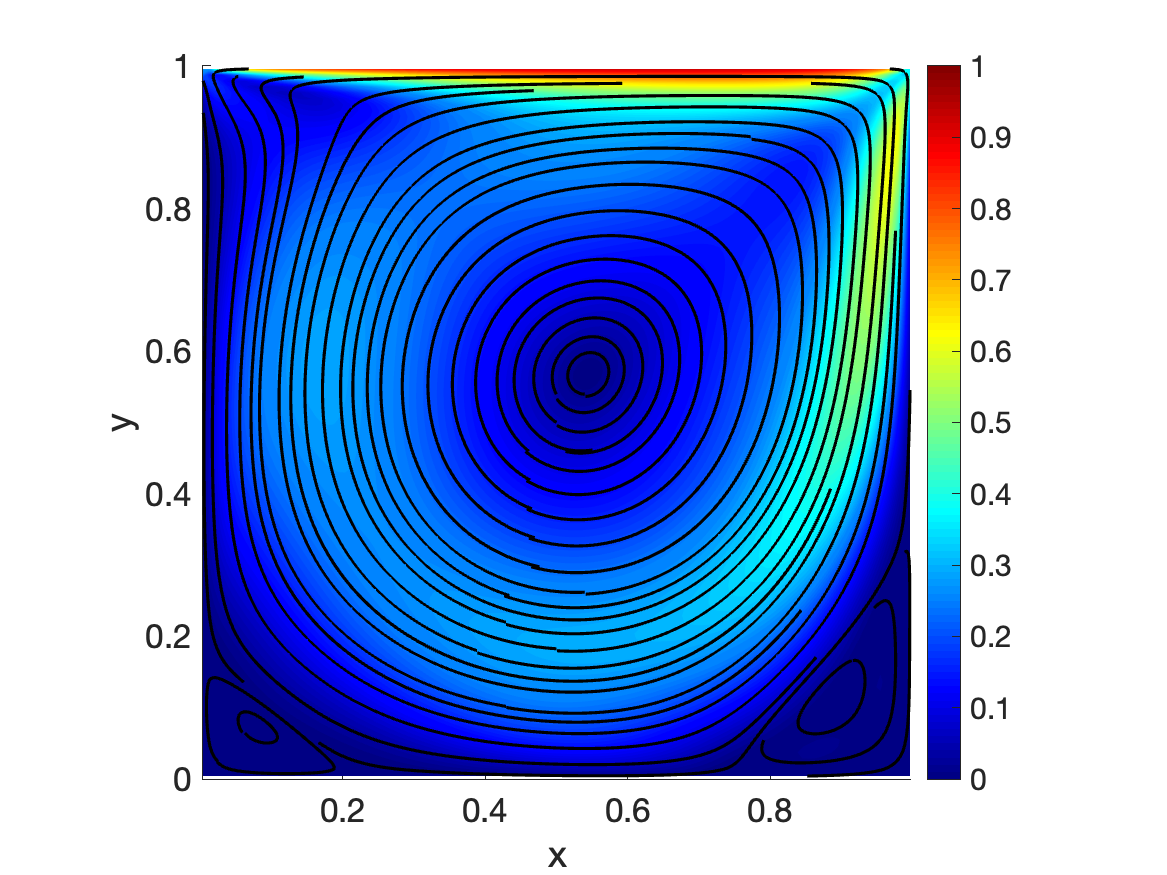}}
\caption{lid-driven cavity: $x$-velocity field and streamline pattern for $\text{Re}=1000$.}
\label{fig:res_cavity_u}
\end{figure}

\subsection{Diffusion problem}
\label{sec:res_diffusion_pb}

We consider a steady-state diffusion problem governed by the equation
\[
-\nabla \cdot (k \nabla T) = Q
\]
in the domain $\Omega = (0, 1) \times (0, 2)$, with the analytical solution
\[
T = 20 + x^2 - xy - 3y^2.
\]
We examine three distinct cases based on varying thermal conductivities, as follows:
\\
\texttt{case 1}: $k = 2T^3-0.1T^2+T$, \\
\texttt{case 2}: $k = -0.1T^2+T$, \\
\texttt{case 3}: $k = 1$. \\
For each case, the source term $ Q $ is computed to ensure that the diffusion equation is satisfied.

\begin{figure}[H]
\centering

\subfloat[exact solution]{\includegraphics[width=0.4\textwidth]{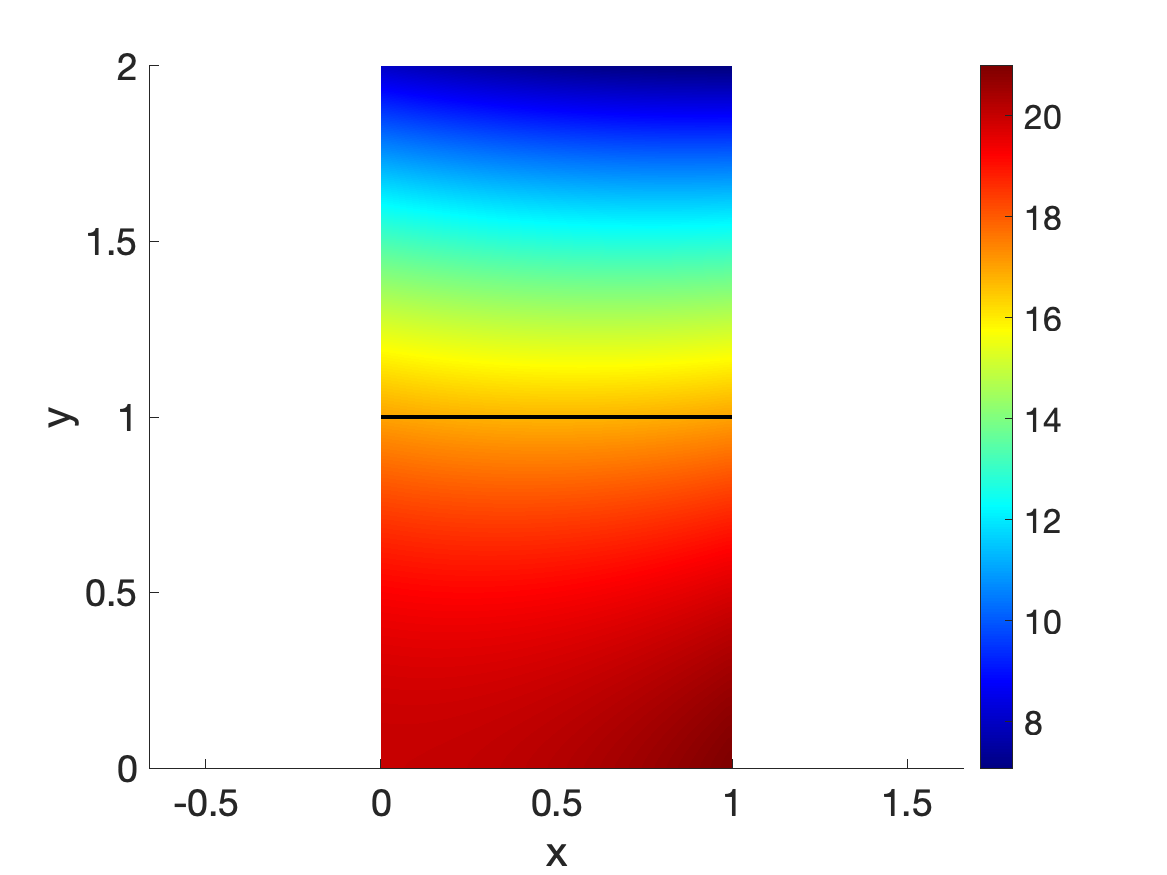}}
~~
\subfloat[error of the DtN method]{\includegraphics[width=0.4\textwidth]{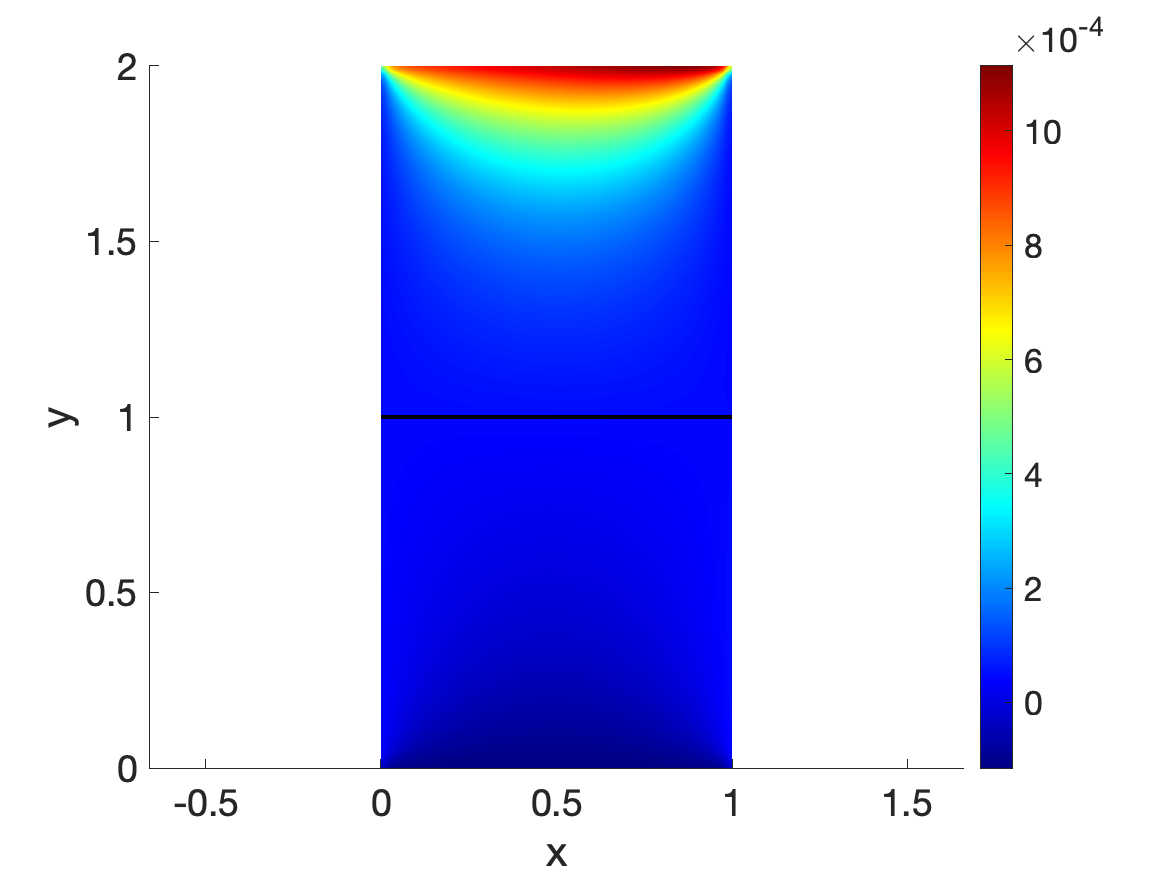}}
\hfill
\subfloat[error of the OB method]{\includegraphics[width=0.4\textwidth]{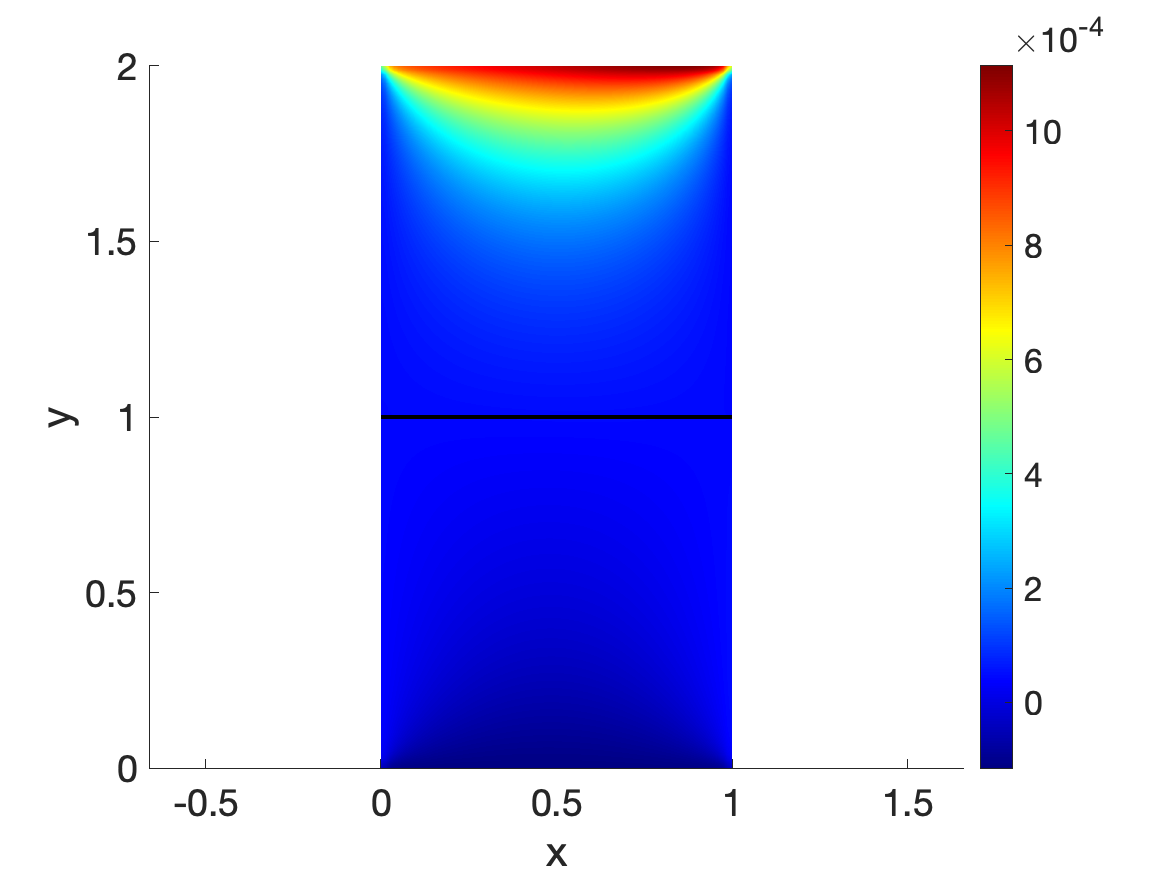}}
~~
\subfloat[difference between the OB method and the reduced OB method]
{\includegraphics[width=0.4\textwidth]{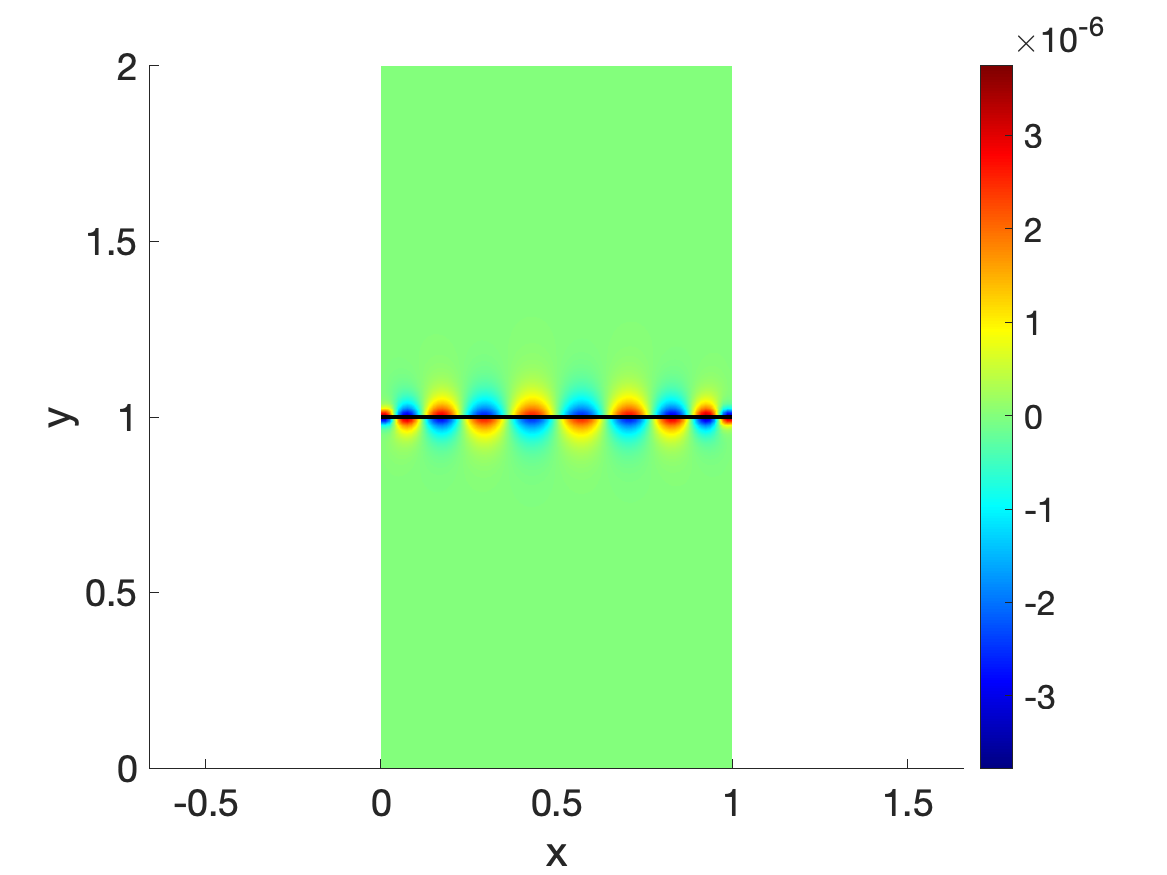}}
\caption{diffusion problem: \texttt{case 2}, $h=\frac{1}{80}$. }
\label{fig:res_diffusion}
\end{figure}

The computational domain is partitioned into two subdomains: a fluid domain and a solid domain, with the interface $\Gamma = (0, 1) \times \{1\}$. In both subdomains, the same diffusion equation is solved using the Dirichlet-to-Neumann (DtN) method and the OB method, with the appropriate coupling conditions as specified in equation \eqref{eq:CHT_coupling_conditions} for the fluid-solid interface. 
Dirichlet boundary conditions are applied to the remaining boundaries. 
The convergence tolerance for both methods is set to $10^{-6}$, and for the DtN method, a relaxation parameter of 0.2 is employed.
For each case, three spatial discretizations are considered, with characteristic cell sizes $h = \frac{1}{20}, \frac{1}{40}, \frac{1}{80}$. The numerical results for \texttt{case 2} with $h = \frac{1}{80}$ are shown in Figure \ref{fig:res_diffusion}, where both the DtN and OB methods yield comparable results, with a maximum pointwise error of $\mathcal{O}(10^{-3})$ relative to the exact solution. 
A line at $y=1$ is drawn to highlight the fluid-solid interface.
Additionally, the test case is computed using the reduced OB method (cf. Section \ref{sec:reduced_OB}) with $N_r = 5$. The difference between the OB method with and without interface reduction is presented in Figure \ref{fig:res_diffusion}(d), where the difference is found to be of $\mathcal{O}(10^{-6})$.

In Table \ref{table:res_diffusion}, we compare the DtN method and the OB method in terms of the number of iterations and computational cost. We remark that the DtN method fails to converge for \texttt{case 1}, even with a relaxation parameter of $0.01$. For all numerical cases, the OB method demonstrates superior performance compared to the DtN method. Additionally, the numerical simulation can be further accelerated using the reduced OB approach, with the acceleration becoming more significant as the spatial discretization is refined. We emphasize that the current code implementation is not optimized, with matrix assembly accounting for a large portion of the computational time. With a more efficient implementation, the reduced OB method should be more effective compared to the OB method. This work serves as a proof of concept, and thus, code optimization is not the primary focus.

\begin{table}[H]
\centering
\begin{tabular}{|c|c|c|c|c|c|c|}
\hline
 & \multicolumn{2}{c|}{{DtN method}} & \multicolumn{2}{c|}{{OB method}} & \multicolumn{2}{c|}{{reduced OB method}} \\ \hline
& number iters & comp cost [s] & number iters & comp cost [s] & number iters & comp cost [s] \\ \hline
{\texttt{case 1} ($h=\frac{1}{20}$) } & - & - & 25 & 0.17 & 25 & 0.16 \\ 
{\texttt{case 1} ($h=\frac{1}{40}$) } & - & - & 26 & 0.94 & 27 & 0.72 \\
{\texttt{case 1} ($h=\frac{1}{80}$) } & - & - & 27 & 5.91 & 27 & 4.43 \\  \hline
{\texttt{case 2} ($h=\frac{1}{20}$) } & 34 & 0.20 & 15 & 0.10 & 15 & 0.08 \\ 
{\texttt{case 2} ($h=\frac{1}{40}$) } & 34 & 0.83 & 16 & 0.50 & 16 & 0.42 \\
{\texttt{case 2} ($h=\frac{1}{80}$) } & 34 & 5.16 & 16 & 3.40 & 16 & 2.62 \\  \hline
{\texttt{case 3} ($h=\frac{1}{20}$) } & 32 & 0.18 & 1 & 0.01 & 1 & 0.01 \\ 
{\texttt{case 3} ($h=\frac{1}{40}$) } & 32 & 0.78 & 1 & 0.06 & 1 & 0.05 \\
{\texttt{case 3} ($h=\frac{1}{80}$) } & 32 & 4.82 & 1 & 0.38 & 1 & 0.32 \\  \hline
\end{tabular}
\caption{comparison of the DtN and the OB methods in terms of number of iterations and computational cost.} 
\label{table:res_diffusion}
\end{table}

\subsection{Flow over a heated plate}

\begin{figure}[H]
\centering
\begin{tikzpicture}

\fill[blue!20] (0,0) rectangle (7,1.5); 
\node at (3.5, 0.75) {$\Omega_{\rm f}$}; 

\fill[red!30] (1.5,-0.8) rectangle (3.5,0); 
\node at (2.5, -0.4) {$\Omega_{\rm s}$}; 

\draw[thick] (0,0) rectangle (7,1.5); 
\draw[thick, red] (1.5,-0.8) rectangle (3.5,0); 

\draw[dashed, thick] (1.5, 0) -- (3.5, 0); 
\node at (2.5, 0.2) {$\Gamma$}; 


\draw[<->, thick] (7.2, 0) -- (7.2, 1.5) node[midway, right] {$0.75$}; 

\draw[<->, thick] (0, -0.7) -- (1.5, -0.7) node[midway, below] {$0.5$}; 

\draw[<->, thick] (3.5, -0.2) -- (7, -0.2) node[midway, below] {$2$}; 

\draw[->, thick] (-1, 1.2) -- (-0.1, 1.2) node[left, xshift=-1cm, yshift=-0.2cm] {$U_{\rm{in}}=0.1$}; 
\draw[->, thick] (-1, 0.3) -- (-0.1, 0.3) node[left, xshift=-1cm, yshift=0.2cm] {$T_{\rm{in}}=300$}; 
\draw[->, thick] (-1, 0.75) -- (-0.1, 0.75) node[left, xshift=-1cm] {}; 

\fill[black] (0,0) circle (2pt); 
\node at (0, 0) [below] {$O(0,0)$}; 

\end{tikzpicture}
\caption{flow over a heated plate: illustration of the computational domain (not scaled). }
\label{fig:illus_flow_over_heated_plate}
\end{figure}
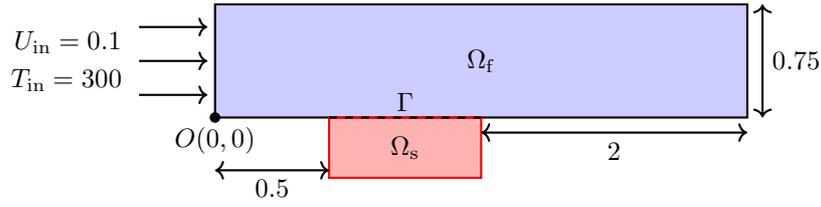

In this test case, we simulate fluid flow over a heated plate \cite{vynnycky1998forced,yau2016} and
investigate the heat transfer characteristics resulting from  convection in the fluid and conduction within the solid plate. The computational domain is depicted in Figure \ref{fig:illus_flow_over_heated_plate}, with the geometric properties specified. The length of the plate, denoted as $\Omega_{\rm s}$, is $L=1$, and the height of the plate is $b=L/4$. 

The fluid enters from the left with a free-stream velocity of $ U_{\text{in}} = 0.1 $. The boundary conditions for the fluid flow are as follows: an outlet condition is specified on the right boundary; a Neumann condition is imposed on the top boundary; a slip wall condition is enforced to the bottom boundary upstream of the leading edge of the plate \cite{vynnycky1998forced}; and a no-slip wall condition is applied to the remaining portion of the bottom boundary.
For heat transfer in the fluid domain, a uniform temperature of $ T_{\text{in}} = 300 $ is imposed at the inlet. All boundaries, except for the interface where the coupling conditions in \eqref{eq:CHT_coupling_conditions} are enforced, are treated as adiabatic, i.e., $ \frac{\partial T}{\partial \mathbf{n}} = 0 $. Regarding heat transfer in the solid plate, the bottom surface is maintained at a uniform temperature of $ T_0 = 310 $, while the left and right sides are assumed to be adiabatic.

The densities of both the fluid and the solid are taken to be equal, with $\rho_{\rm f} = \rho_{\rm s} = 1$. The dynamic viscosity of the fluid is chosen to be $\mu_{\rm f} = 2 \times 10^{-4}$, which results in a Reynolds number $\text{Re} = \frac{\rho_{\rm f} U_{\rm{in}} L}{\mu_{\rm f}} = 500$, based on the characteristic length of the plate $L$. The thermal properties of the solid plate include a thermal conductivity of $k_{\rm s} = 100$ and a specific heat capacity of $c_{\rm{ps}} = 100$.
For the fluid, the thermal conductivity $k_{\rm f}$ and the specific heat capacity $c_{\rm{pf}}$ are computed using the conductivity ratio $k = \frac{k_{\rm s}}{k_{\rm f}}$ and the Prandtl number $\text{Pr} = \frac{\mu_{\rm f} c_{\rm{pf}}}{k_{\rm f}}$. The effects of these parameters are investigated by varying the conductivity ratio, $k = 1, 2, 5, 20$, and the Prandtl number, $\text{Pr} = 0.01, 100$.

Note that gravity is not considered in this case, resulting in the decoupling of the energy equation from the mass and momentum equations of the fluid. As a consequence, the fluid flow $\mathbf{u}_{\rm f}$ can be determined by solving a purely fluid dynamics problem \eqref{eq:fluid_mass}-\eqref{eq:fluid_momentum}. Subsequently, the heat transfer problem can be addressed by solving the governing equations \eqref{eq:CHT_fluid} and \eqref{eq:CHT_solid}, with $\mathbf{u}_{\rm f}$ already known.

\begin{figure}[H]
\centering

\subfloat[semi-implicit method]{\includegraphics[width=0.45\textwidth]{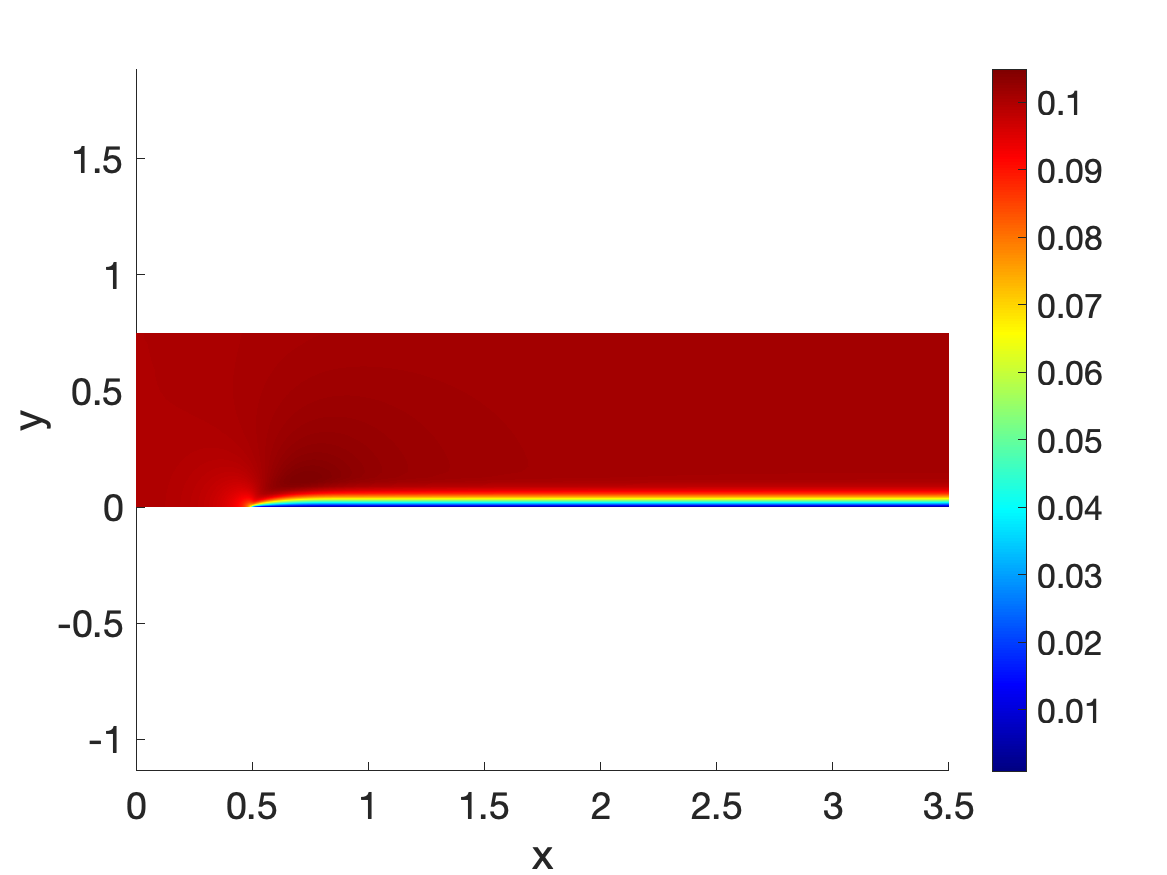}}
~~
\subfloat[difference between the semi-implicit method and the SIMPLE method]{\includegraphics[width=0.45\textwidth]{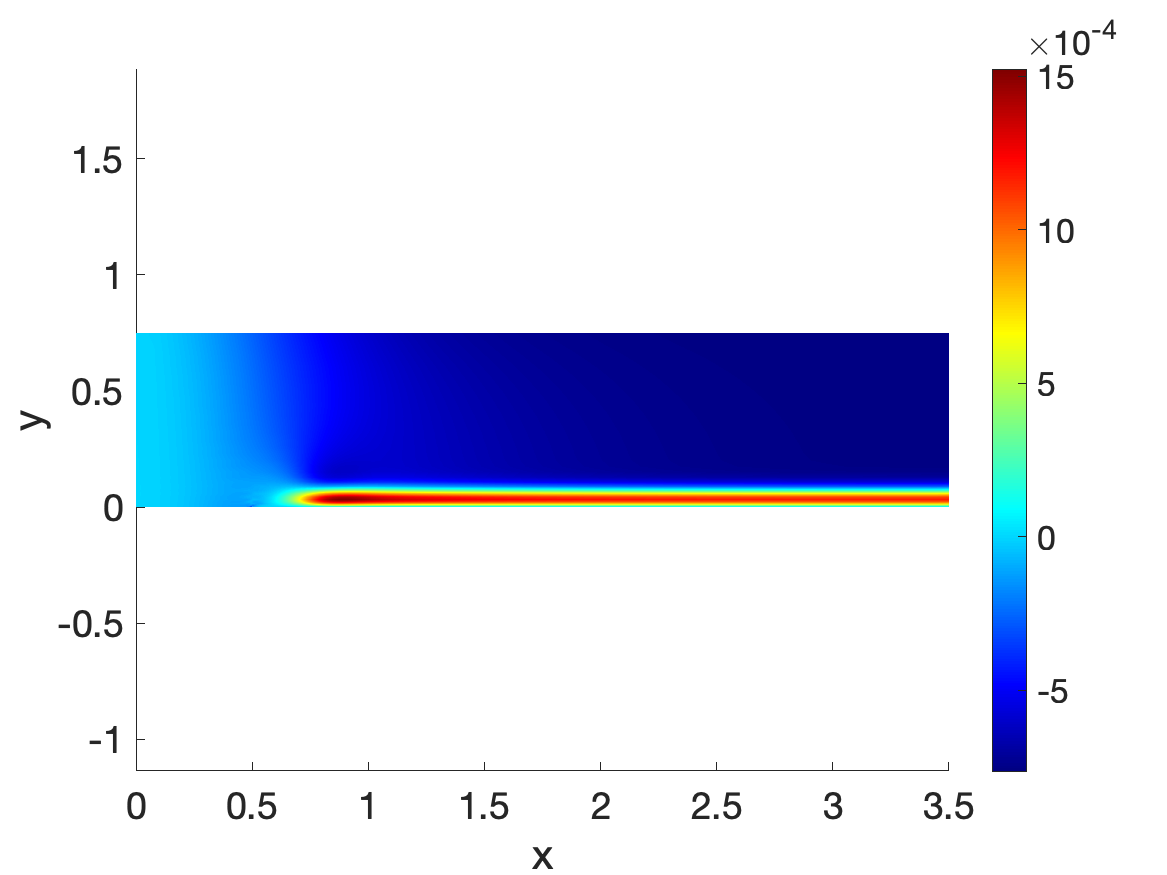}}
\caption{flow over a heated plate: $x$-velocity field at steady state. }
\label{fig:res_flow_plate_u}
\end{figure}

\begin{figure}[H]
\centering
{\includegraphics[width=8cm]{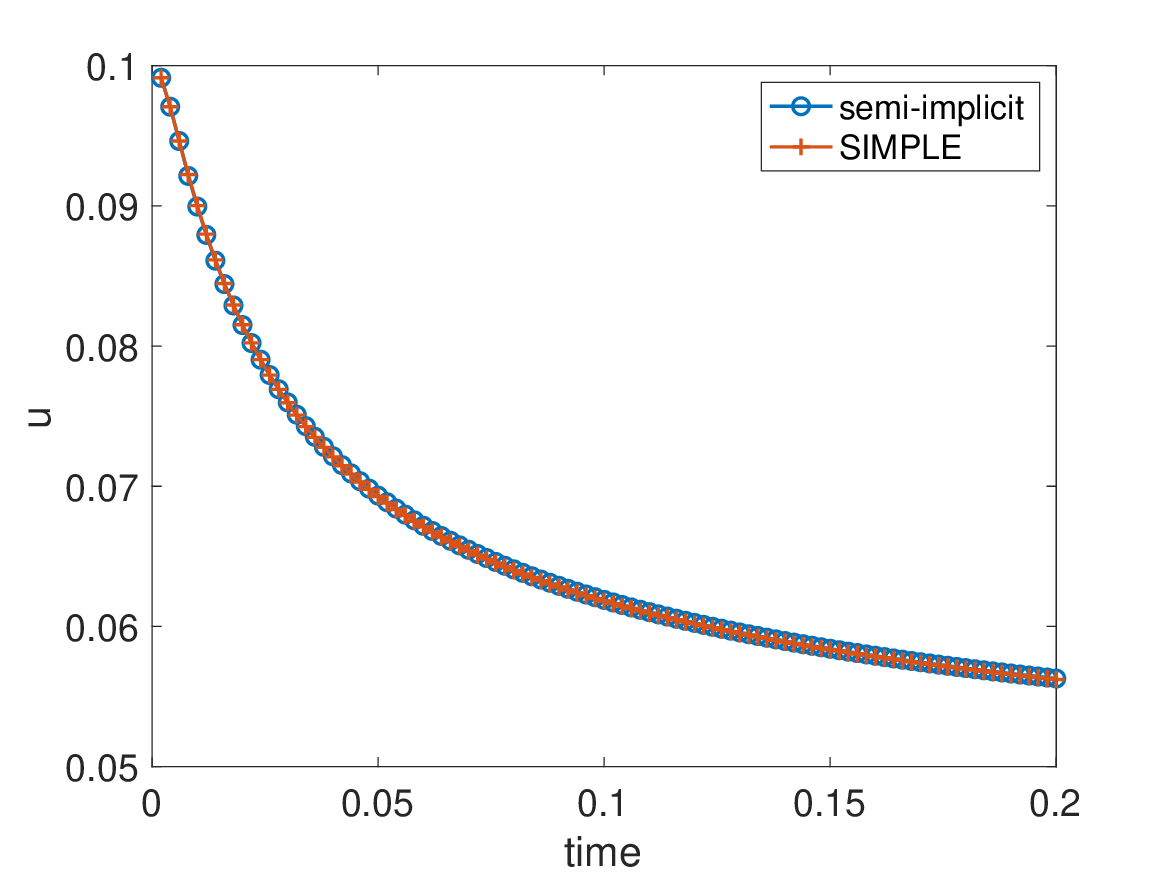}}
\caption{flow over a heated plate: evolution of the velocity at a probe point.}
\label{fig:res_flow_plate_u_probe}
\end{figure}

The fluid flow is solved using both the SIMPLE method and the proposed semi-implicit method. The resulting velocity fields in the $x$-direction  are presented for both methods at steady-state conditions in Figure \ref{fig:res_flow_plate_u}. The steady state is considered to be achieved when the relative change in velocity between two successive time steps is less than $10^{-5}$ (cf. equation \eqref{eq:steady_state}).
We compare the evolution of $x$-direction velocity at a probe point $(0.5, 0.002)$, located near the leading edge of the solid plate, obtained using both the SIMPLE and semi-implicit methods in Figure \ref{fig:res_flow_plate_u_probe}. The results from both methods exhibit good agreement. The SIMPLE method requires $4-5$ iterations per time step during the initial stages and $2-3$ iterations when approaching steady state. In contrast, the semi-implicit method, which does not require iterations (i.e., we set $K = 0$ for the iterative solver described in \ref{appendix:iter_solver}), demonstrates higher efficiency.

\begin{figure}[H]
    \centering
    \begin{minipage}{0.49\textwidth}
        \centering
        \includegraphics[width=\textwidth]{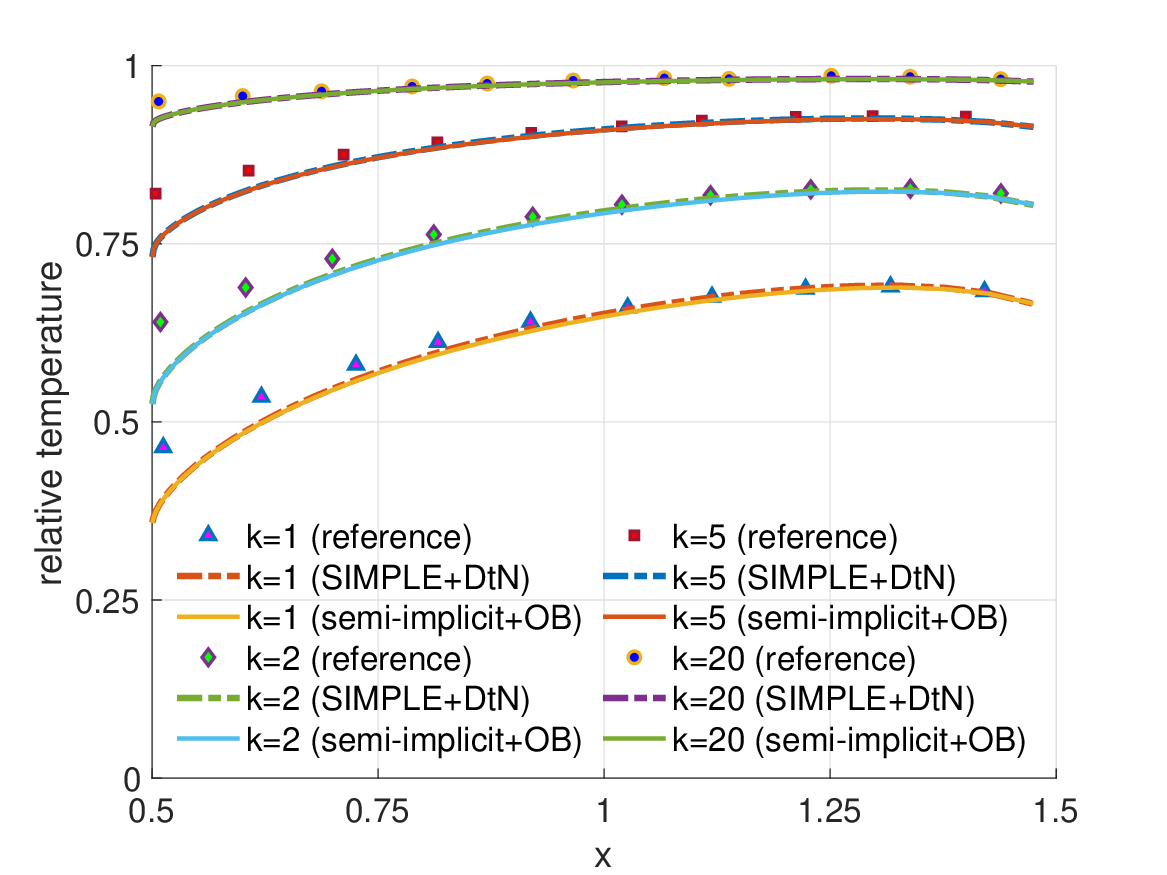} 
        \caption*{$\text{Pr}=0.01$}
        \label{fig:image1}
    \end{minipage}
    \hfill
    \begin{minipage}{0.49\textwidth}
        \centering
        \includegraphics[width=\textwidth]{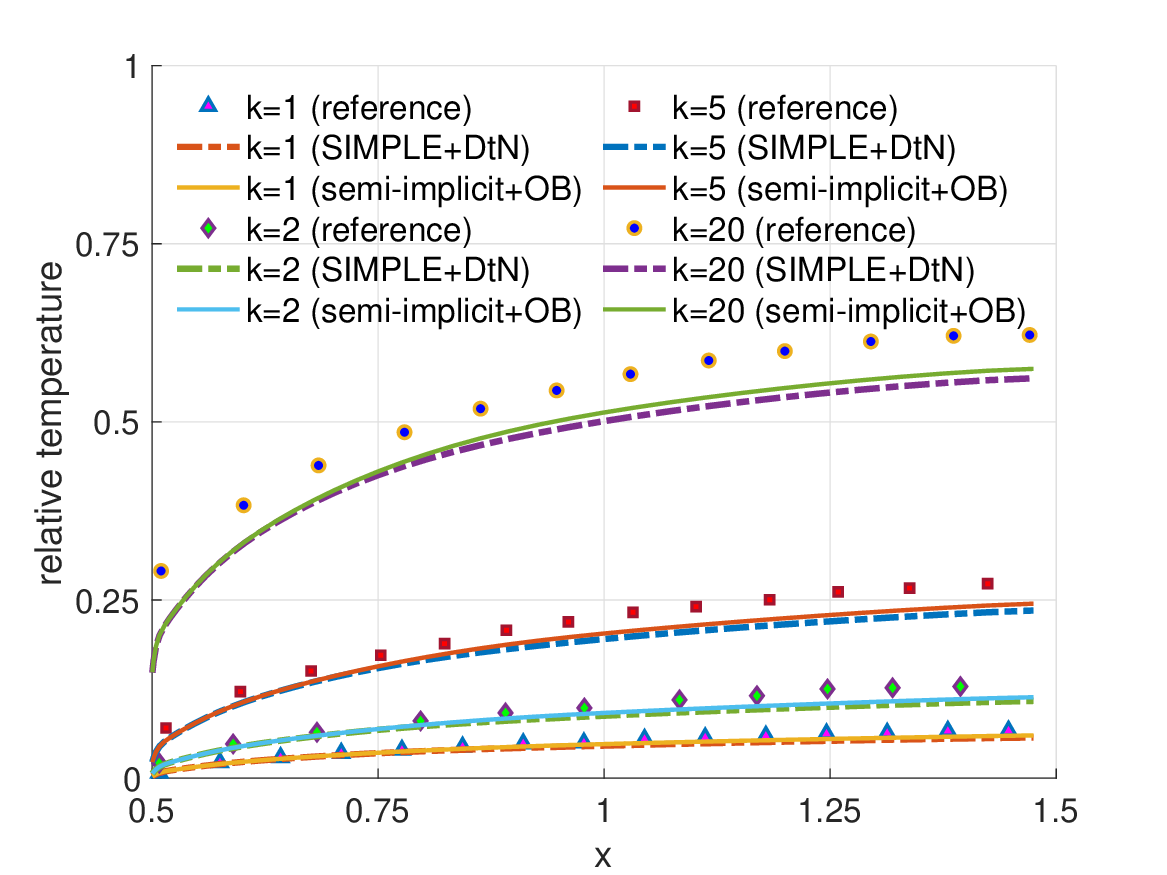}
        \caption*{$\text{Pr}=100$}
        \label{fig:image2}
    \end{minipage}
    \caption{flow over a heated plate: relative temperature distribution along the fluid-solid interface for various values of $k$ and $\text{Pr}$.}
    \label{fig:flow_plate_compr_T}
\end{figure}

Once the fluid flow is solved using either the semi-implicit method or the SIMPLE method, the heat transfer between the fluid and the solid plate can subsequently be addressed using either the DtN method or the OB method. In Figure \ref{fig:flow_plate_compr_T}, we present the relative temperature distribution along the fluid-solid interface for both the DtN and OB methods, compared against the reference results provided in \cite{vynnycky1998forced} for various conductivity ratios $k$ and Prandtl numbers $\text{Pr}$, where 
the relative temperature is defined as:
\[
T_{\text{rel}} = \frac{T - T_{\text{in}}}{T_0 - T_{\text{in}}}.
\]
Both methods yield results in close agreement with each other and the reference values. In this test case, as the heat transfer model is linear, the OB method converges within a single iteration, similar to \texttt{case 3} in the diffusion problem discussed in Section \ref{sec:res_diffusion_pb}. For the case with $\text{Pr} = 0.01$, we use a relaxation parameter of $0.2$ for the DtN method, which requires  approximately 35--65 iterations to attain convergence, for different choices of $k$. The OB method achieves a speed-up factor of 12--22 compared to the DtN method, highlighting its superior efficiency. 
For the case with $\text{Pr} = 100$, the efficiency of the OB method is even more pronounced, as the DtN method exhibits convergence difficulties, requiring a significantly reduced relaxation parameter of $0.05$ to ensure convergence.

\subsection{Natural convection}

We consider natural convection in a square enclosure $\Omega_{\rm f}$, as illustrated in Figure \ref{fig:CHT_illustration}, with a length of $L = 1$. The right side is heated by a solid domain $\Omega_{\rm s}$ of width $b = 0.2$. Since the buoyancy term is included, Algorithm \ref{algo:OB_CHT} is applied to solve this problem. 
For the fluid domain, a no-slip wall condition is imposed for the fluid flow. The top and bottom boundaries are considered adiabatic, while the left boundary is maintained at a cold temperature $T_{\rm c} = 1$. In the solid domain, the top and bottom boundaries are also adiabatic, whereas the right boundary is set to a hot temperature $T_{\rm h} = 2$. The gravitational acceleration is given by $\mathbf{g}_0 = (0,\,-1)^{\top}$, and the reference temperature in the buoyancy term is chosen as $T_{\rm{ref}} = T_{\rm c}$.

We consider three parameter settings \cite{kazemi2013high,kazemi2014analysis,pan2018efficient}: \\
\texttt{case 1}: $\rho_{\rm f} = 1$, $\mu_{\rm f} = 0.7$, $\beta = 0.7 \times 10^5$, $\rho_{\rm s} = 7.5 \times 10^3$, $k_{\rm s} = 1.6 \times 10^3$, $c_{\rm{ps}} = 0.5$, $k = 1.6 \times 10^3$, ${\rm {Pr}} = 0.7$. \\ 
\texttt{case 2}: $\rho_{\rm f} = 1$, $\mu_{\rm f} = 7$, $\beta = 4.9 \times 10^5$, $\rho_{\rm s} = 7.5$, $k_{\rm s} = 80$, $c_{\rm{ps}} = 0.12$, $k = 80$, ${\rm {Pr}} = 7$. \\ 
\texttt{case 3}: $\rho_{\rm f} = 1$, $\mu_{\rm f} = 7$, $\beta = 4.9 \times 10^5$, $\rho_{\rm s} = 7.5$, $k_{\rm s} = 2.7$, $c_{\rm{ps}} = 0.0576$, $k = 2.7$, ${\rm {Pr}} = 7$. \\ 
The Prandtl number is given by
$
\text{Pr} = \frac{\mu_{\rm f} c_{\rm{pf}}}{k_{\rm f}}.
$
These three test cases correspond to an increasing ratio of the thermal effusivities of the coupled domains, defined as
\[
\sigma = \frac{k_{\rm f}}{k_{\rm s}} \sqrt{\frac{\alpha_{\rm s}}{\alpha_{\rm f}}},
\]
where the thermal diffusivity is defined as
$
\alpha_{(\cdot)} = \frac{k_{(\cdot)}}{\rho_{(\cdot)} c_{\rm p(\cdot)}}, \; (\cdot) = \rm{f, s}.
$
The thermal effusivity ratios are $\sigma = 3.95 \times 10^{-4}$, $0.12$, and $0.93$, respectively,  indicating weak thermal interaction between the fluid and solid domains in \texttt{case 1}, strong interaction in \texttt{case 3}, and an intermediate level of interaction in \texttt{case 2}.

For this test case, we rely on the iterative semi-implicit method, as outlined in Algorithm \ref{alg:semi_implicit_iter}, to solve the fluid flow. We use a spatial discretization of size $\Delta x=\Delta y = \frac{1}{80}$, and a time step of $\Delta t=0.001$. The number of iterations is slightly smaller  than that required by the SIMPLE method. However, a key advantage is that the same matrix is solved  at each iteration and time step. 
The outer loop of the optimization-based algorithm (cf. Algorithm \ref{algo:OB_CHT}) typically requires $5$ iterations to converge in the initial stages and $2$ iterations as it approaches the steady state. Since the heat transfer problem is linear in this test case, the inner optimization problem converges within $1$ sub-iteration. In contrast, if we use the DtN method instead of the OB method in step $4$ of Algorithm \ref{algo:OB_CHT}, the DtN method typically requires $20-40$ sub-iterations. These results demonstrate the improved efficiency of the proposed approach.

\begin{figure}[H]
\centering

\subfloat[$u$]{\includegraphics[width=0.32\textwidth]{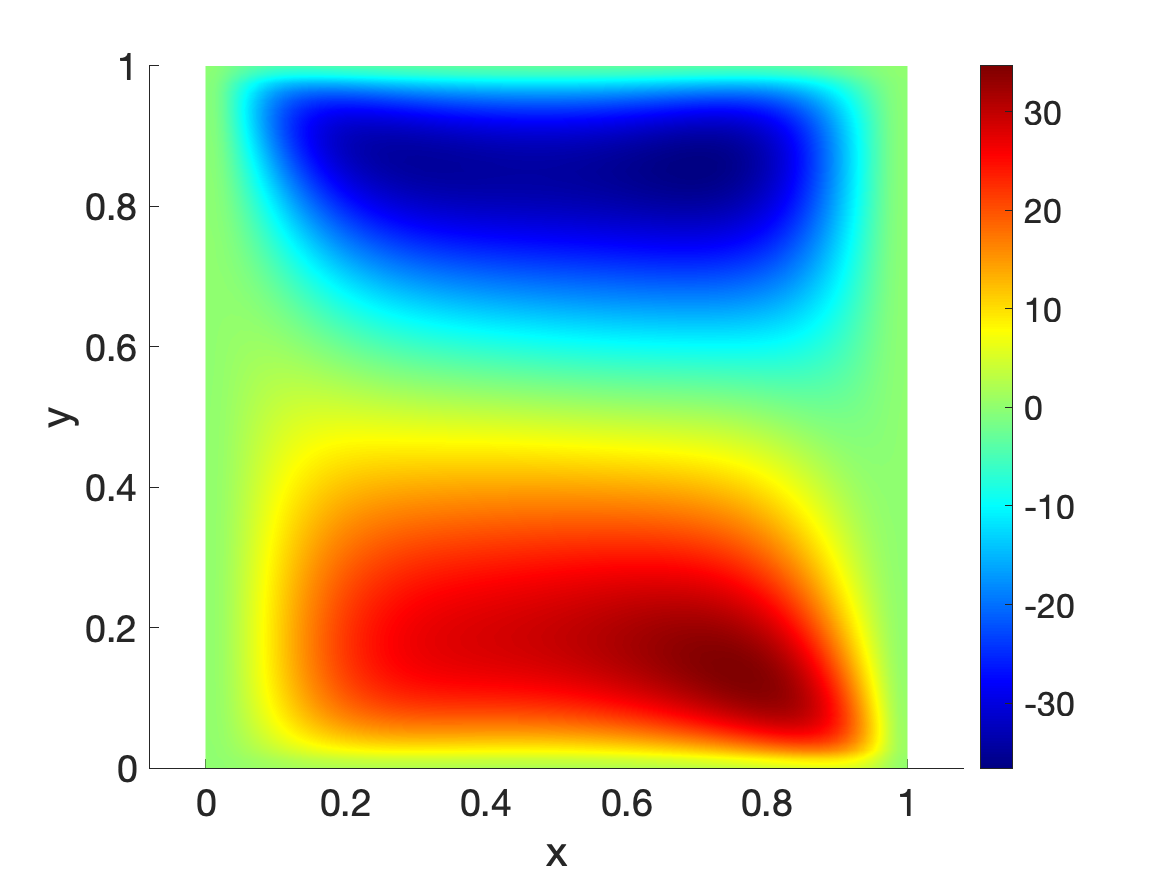}}
~~
\subfloat[$v$]{\includegraphics[width=0.32\textwidth]{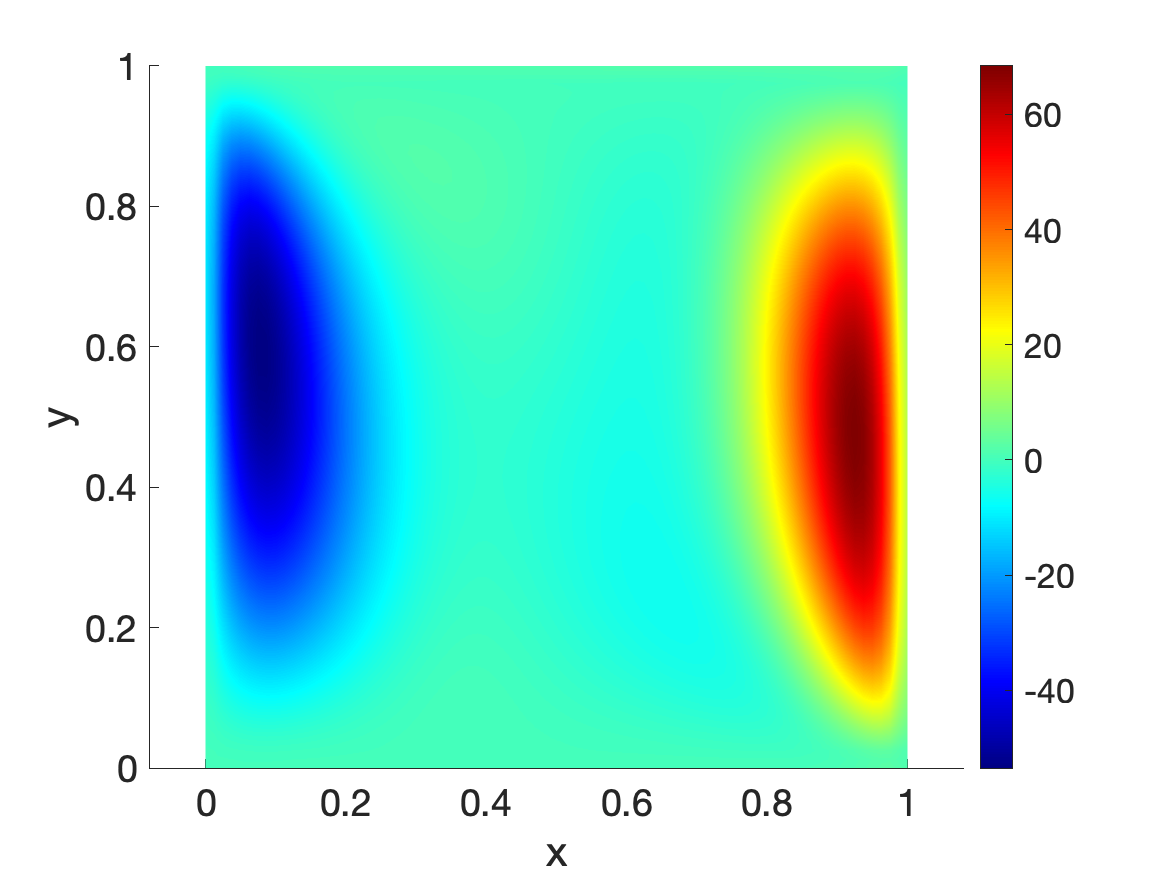}}
~~
\subfloat[$T$]{\includegraphics[width=0.32\textwidth]{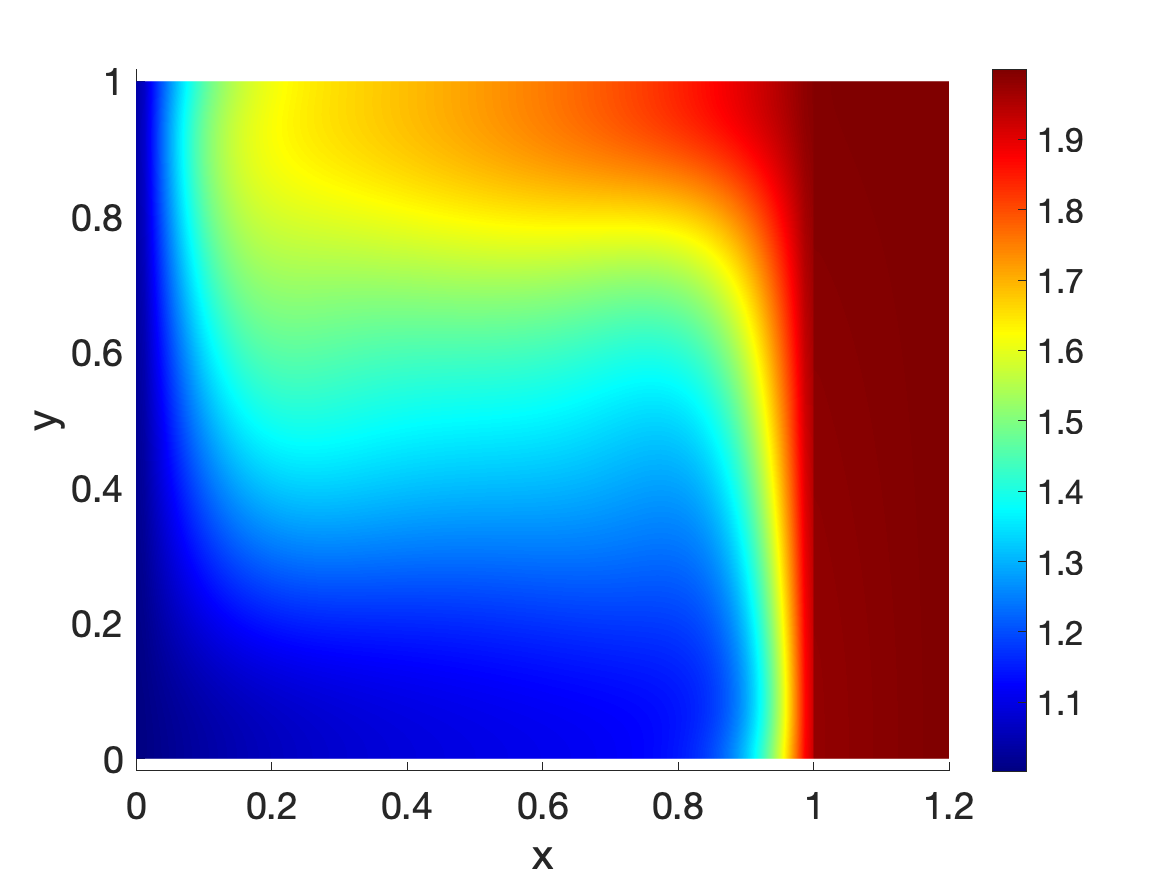}}
\caption{natural convection problem: \texttt{case 2}, solution fields at $t=0.07$. }
\label{fig:natural_conv_fields}
\end{figure} 

\begin{figure}[H]
\centering

\subfloat[temperature profile along $y=0.5$ at $t=0.07$]{\includegraphics[width=0.42\textwidth]{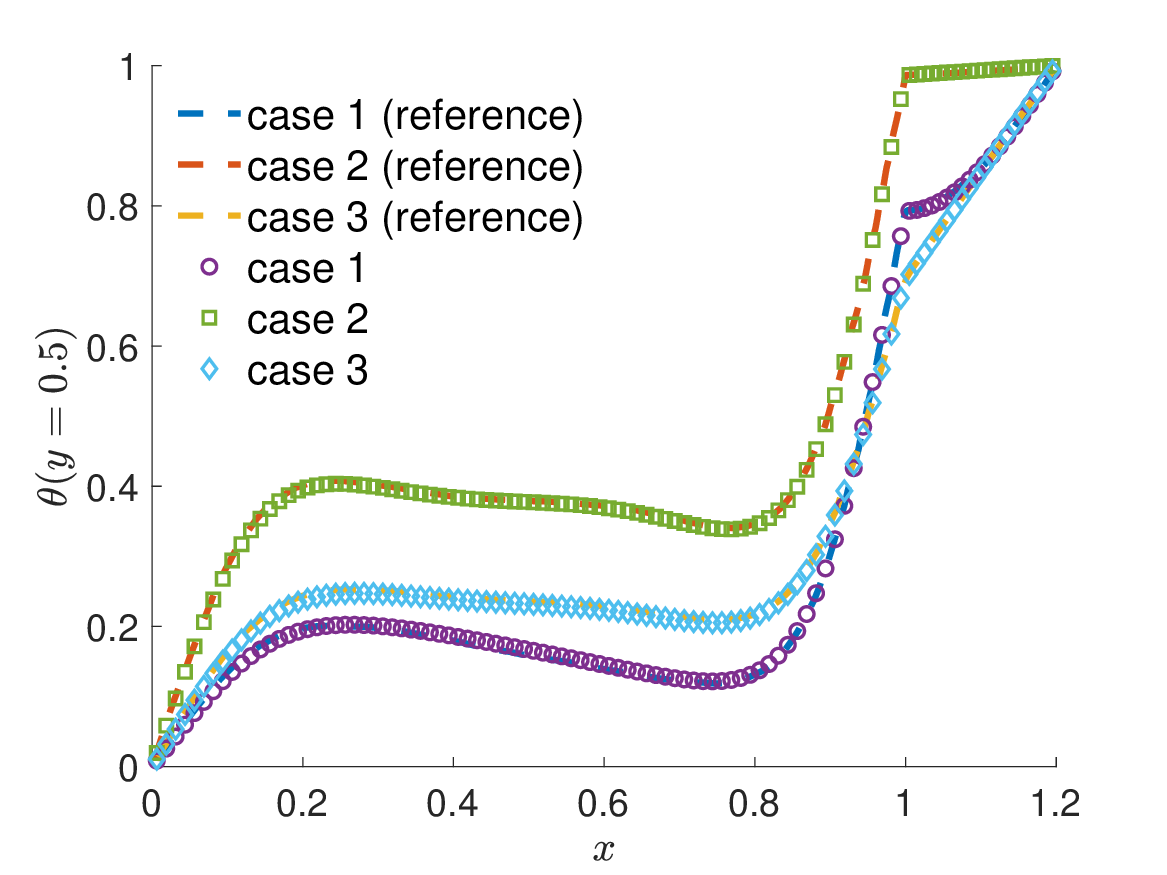}}
~~
\subfloat[temperature evolution at the center of the fluid domain]{\includegraphics[width=0.42\textwidth]{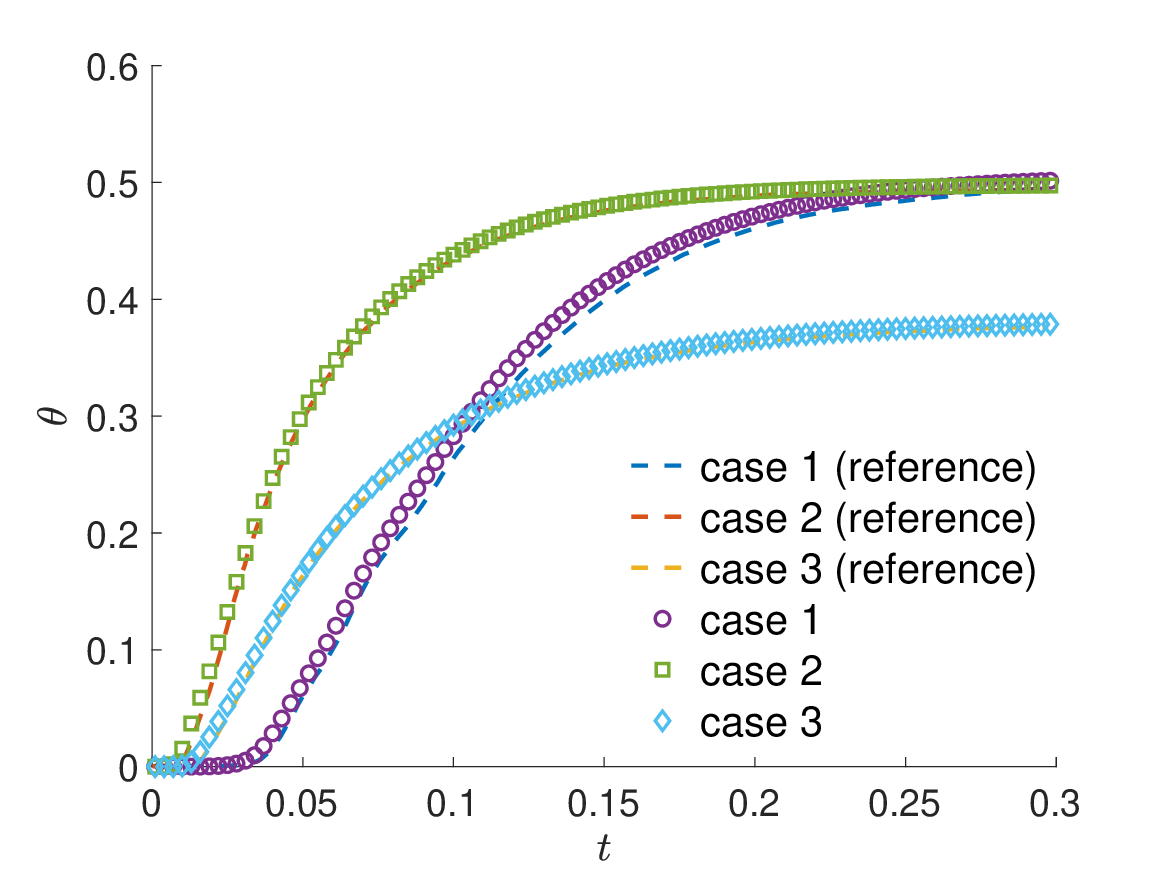}}
\caption{natural convection problem: temperature profile. }
\label{fig:natural_conv_profile_T}
\end{figure}

We present the solution fields at $t = 0.07$ for \texttt{case 2} in Figure \ref{fig:natural_conv_fields}. As observed, the fluid circulates in the counter-clockwise direction. The flow moves upward along the hot fluid-solid interface, leading to an increase in temperature in the $y$-direction at the interface. Similarly, the flow moves downward along the cold left boundary, causing a temperature decrease in the negative $y$-direction at the left boundary. Similar behavior is observed in the other two cases.
We present the dimensionless temperature profile $\theta$ in Figure \ref{fig:natural_conv_profile_T} for the three test cases, defined as
\begin{align*}
\theta = \frac{T - T_{\rm c}}{T_{\rm h} - T_{\rm c}}.
\end{align*}
Specifically, Figure \ref{fig:natural_conv_profile_T}(a) shows the temperature profile along $y = 0.5$ at $t = 0.07$. The obtained results exhibit good agreement with the reference data \cite{pan2018efficient}. In Figure \ref{fig:natural_conv_profile_T}(b), we present the time evolution of the temperature at the center of the fluid domain, located at $(0.5, 0.5)$. The proposed approach accurately captures the dynamic behavior.

\section{Conclusions}
We proposed a partitioned numerical framework for solving conjugate heat transfer problems by integrating a semi-implicit method for fluid flow with an optimization-based approach for heat transfer across the fluid-solid interface. The semi-implicit method efficiently solves the incompressible Navier-Stokes equations while maintaining stability under a Courant-Friedrichs-Lewy (CFL) condition, with stability condition rigorously established through von Neumann analysis. The optimization-based method effectively enforces fluid-solid interface conditions and efficiently solves heat transfer between the fluid and solid domains. The framework improves computational performance by eliminating sub-iterations within each time step for solving the Navier-Stokes equations and by leveraging a sequential quadratic programming approach to solve the constrained optimization problem. Numerical results validate the proposed method, showing improved efficiency compared to the SIMPLE method for fluid flow and the Dirichlet-to-Neumann method for heat transfer. 

In the future, we plan to develop efficient factorization techniques for solving the semi-implicit system to enhance computational efficiency;   we plan to extended the method towards a nearly implicit scheme \cite{trapp1986nearly}, enabling for large time steps and steady-state computations. Additionally, the optimization-based method will be applied to more challenging problems, such as those involving complex geometries or highly nonlinear thermal interactions, to assess its robustness in real-world engineering applications. Furthermore, the introduction of Robin-type boundary conditions in the optimization-based formulation will be explored  to improve interface coupling and accelerate convergence in partitioned solvers.

\section*{Acknowledgements}
The work of Xiandong Liu and Lei Zhang is supported by the Fundamental Research Funds for the Central Universities of Tongji University, and National Natural Science Foundation of China (grant number: 12301557). The work of Liang Fang and Lei Zhang is supported by the Interdisciplinary Joint Research Project of Tongji University.

\appendix

\section{Numerical stability}
\label{appendix:stability}

We derive the stability condition of the proposed semi-implicit method presented in \ref{sec:semi_implicit}. Following the approach in \cite{chorin1997numerical}, we consider the non-dimensional form of the governing equations and introduce an  artificial compressibility term into the mass conservation equation. Additionally, to facilitate analysis,  we include a convective term in the mass conservation equation, thereby reformulating the mass conservation equation into a form analogous to that of compressible flows. 
The stability properties are examined using von Neumann stability analysis applied to a two-dimensional problem on structured grids. 

With the non-dimensional formulation, the introduction of the artificial compressibility term, and the inclusion of the convective term, the Navier-Stokes equations considered in this analysis are given as
\begin{align}
\label{eq:mass_ns_stab}
    \frac{\partial \rho}{\partial t}
    +\rho\frac{\partial  u}{\partial x}+\rho\frac{\partial  v}{\partial y}
+u\frac{\partial \rho }{\partial x}+v\frac{\partial \rho }{\partial y}    
    &=0, \\
        \frac{\partial u}{\partial t} + u\frac{\partial u}{\partial x}+v\frac{\partial u}{\partial y}-\mu\left(\frac{\partial^2 u}{\partial x^2}+\frac{\partial^2 u}{\partial y^2}\right)+\frac{\partial p}{\partial x}&=0, \\
            \frac{\partial v}{\partial t} + v\frac{\partial v}{\partial y}+u\frac{\partial v}{\partial x}-\mu\left(\frac{\partial^2 v}{\partial x^2}+\frac{\partial^2 v}{\partial y^2}\right)+\frac{\partial p}{\partial y}&=0,
\label{eq:momentum_ns_stab}
\end{align}
where the artificial density $\rho$ is related to the dimensionless pressure through the artificial speed of sound $c$ as
\begin{equation*}
    \rho=\frac{p}{c^2}.
\end{equation*}
Assuming the velocity components satisfy $u > 0$ and $v > 0$, we apply the semi-implicit to the above system \eqref{eq:mass_ns_stab}-\eqref{eq:momentum_ns_stab}. By linearizing the system  around a reference state $(p_0, u_0, v_0, \rho_0, c_0)$, we obtain the following discretized equations 
\begin{align}\label{eq:mass_linearized}
&    p_{j,k}^{n+1}-p_{j,k}^{n}
    +\frac{\rho_0c_0^2\Delta t}{\Delta x}\left(u^{n+1}_{j+\frac{1}{2},k}-u^{n+1}_{j-\frac{1}{2},k}\right)
    +\frac{\rho_0c_0^2\Delta t}{\Delta y}\left(v^{n+1}_{j,k+\frac{1}{2}}-v^{n+1}_{j,k-\frac{1}{2}}\right) \nonumber \\
& \quad\quad\quad\quad\quad\quad
\quad\quad\quad\quad\quad\quad  +\frac{u_0\Delta t}{\Delta x}\left(p^n_{j,k}-p^n_{j-1,k}\right)
    +\frac{v_0\Delta t}{\Delta y}\left(p^n_{j,k}-p^n_{j,k-1}\right)    
    =0, \\    
&    u_{j,k}^{n+1}-u_{j,k}^{n}+\frac{u_{0}\Delta t}{\Delta x}\left(u^{n}_{j,k}-u^{n}_{j-1,k}\right)+\frac{v_{0}\Delta t}{\Delta y}\left(u^{n}_{j,k}-u^{n}_{j,k-1}\right)\nonumber \\[5pt]
& \quad\quad\quad    -\frac{2\mu\Delta t}{(\Delta x)^2}\frac{u^{n+1}_{j+1,k}-2u^{n+1}_{j,k}+u^{n+1}_{j-1,k}}{2}-\frac{2\mu\Delta t}{(\Delta y)^2}\frac{u^{n+1}_{j,k+1}-2u^{n+1}_{j,k}+u^{n+1}_{j,k-1}}{2}+\frac{\Delta t}{2\Delta x}\left(p^{n+1}_{j+1,k}-p^{n+1}_{j-1,k}\right)=0, \\
&    v_{j,k}^{n+1}-v_{j,k}^{n}+\frac{v_{0}\Delta t}{\Delta y}\left(v^{n}_{j,k}-v^{n}_{j,k-1}\right)+\frac{u_{0}\Delta t}{\Delta x}\left(v^{n}_{j,k}-v^{n}_{j-1,k}\right)\nonumber \\[5pt]
& \quad\quad\quad    -\frac{2\mu\Delta t}{(\Delta x)^2}\frac{v^{n+1}_{j+1,k}-2v^{n+1}_{j,k}+v^{n+1}_{j-1,k}}{2}-\frac{2\mu\Delta t}{(\Delta y)^2}\frac{v^{n+1}_{j,k+1}-2v^{n+1}_{j,k}+v^{n+1}_{j,k-1}}{2}+\frac{\Delta t}{2\Delta y}\left(p^{n+1}_{j,k+1}-p^{n+1}_{j,k-1}\right)=0.
\end{align}
Moreover, the Rhie-Chow interpolation is employed to compute the velocity at cell face in the mass equation \eqref{eq:mass_linearized}
\begin{align}\label{RC_2D}
&u_{j+\frac{1}{2},k}^{n+1}=\frac{1}{2}(u_{j,k}^{n+1}+u_{j+1,k}^{n+1})-\frac{\Delta t}{(1+C_\mu)}\frac{p_{j+1,k}^{n+1}-p_{j,k}^{n+1}}{\Delta x} +\frac{\Delta t}{2(1+C_\mu)}\left(\frac{p_{j+1,k}^{n+1}-p_{j-1,k}^{n+1}}{2\Delta x} +\frac{p_{j+2,k}^{n+1}-p_{j,k}^{n+1}}{2\Delta x}\right), \\
&v_{j,k+\frac{1}{2}}^{n+1}=\frac{1}{2}(v_{j,k}^{n+1}+v_{j,k+1}^{n+1})-\frac{\Delta t}{(1+C_\mu)}\frac{p_{j,k+1}^{n+1}-p_{j,k}^{n+1}}{\Delta y} +\frac{\Delta t}{2(1+C_\mu)}\left(\frac{p_{j,k+1}^{n+1}-p_{j,k-1}^{n+1}}{2\Delta y} +\frac{p_{j,k+2}^{n+1}-p_{j,k}^{n+1}}{2\Delta y}\right),
\end{align}
where $C_{\mu}=\frac{2\mu\Delta t}{(\Delta x)^2}+\frac{2\mu\Delta t}{(\Delta y)^2}$. 

A von Neumann analysis can be performed by substituting $ \phi^n_{j,k} $ with $ \phi e^{-ij\theta_1}e^{-ik\theta_2} g^n$, 
where $g $ represents the eigenvalue of the amplification matrix. This leads to the following system of equations
\begin{align*}
&    p\left[g-1+C_{m1}k_{11}+C_{m2}k_{12}+\left(\rho_0C_{s1}^2C_{\mu3}k_{31}+\rho_0C_{s2}^2C_{\mu3}k_{32}\right)g \right]
    +u\big(\rho_0c_0C_{s1}k_{21} \big)g
    +v\big(\rho_0c_0C_{s2}k_{22} \big)g=0, \\
&    p\left(c_0^{-1}C_{s1}k_{21}\right)g
+u\left[g-1 +C_{m1}k_{11}+C_{m2}k_{12}-\left(C_{\mu1}k_{41}+C_{\mu2}k_{42}\right)g\right]=0, \\
&    p\left(c_0^{-1}C_{s2}k_{22}\right)g
+v\left[g-1 +C_{m2}k_{12}+C_{m1}k_{11}-\left(C_{\mu1}k_{41}+C_{\mu2}k_{42}\right)g\right]=0,
\end{align*}
where
\begin{align*}
&    C_{m1}=\frac{u_0\Delta t}{\Delta x},\; C_{m2}=\frac{v_0\Delta t}{\Delta y},\; C_{s1}=\frac{c_0\Delta t}{\Delta x},\; C_{s2}=\frac{c_0\Delta t}{\Delta y},\; C_{\mu1}=\frac{2\mu\Delta t}{(\Delta x)^2},\; C_{\mu2}=\frac{2\mu\Delta t}{(\Delta y)^2},\; C_{\mu3}=\frac{1}{1+C_{\mu}}, \\
&    k_{1l}=1-e^{i\theta_l}=1-\cos{\theta_l}-i\sin{\theta_l},\quad k_{2l}=\frac{e^{-i\theta_l}-e^{i\theta_l}}{2}=-i\sin{\theta_l},\\
&    k_{3l}=\frac{e^{2i\theta_l}-4e^{i\theta_l}+6-4e^{-i\theta_l}+e^{-2i\theta_l}}{4}=(\cos{\theta_l}-1)^2,\quad k_{4l}=\frac{e^{-i\theta_l}-2+e^{i\theta_l}}{2}=\cos{\theta_l}-1,\quad l=1,2. 
\end{align*}
In matrix form, this system can be written as:
\begin{align*}
\begin{pmatrix}
G_{11}(g) & g\rho_0c_0C_{s1}k_{21} &   g\rho_0c_0C_{s2}k_{22} \\
gc_0^{-1}C_{s1}k_{21}  & G_{22}(g) &  0\\
gc_0^{-1}C_{s2}k_{22} &  0  &   G_{22}(g)
\end{pmatrix} 
\begin{pmatrix}
p \\
u \\
v
\end{pmatrix}=
\begin{pmatrix}
0 \\
0 \\
0
\end{pmatrix},
\end{align*}
where 
\begin{align*}
G_{11}(g) &= g(1 + \rho_0 C_{s1}^2 C_{\mu3} k_{31} +
\rho_0 C_{s2}^2 C_{\mu3} k_{32}) + G_k, \\
G_{22}(g) &= g(1 - C_{\mu1} k_{41} - C_{\mu2} k_{42}) + G_k, \\
G_k &= C_{m1} k_{11} + C_{m2} k_{12} - 1.
\end{align*}
The characteristic polynomial can be computed as:
\begin{equation}
    \phi_3(g)=\phi_1(g)\phi_2(g), 
    \label{eq:charac_poly}
\end{equation}
where the subscripts of $\phi$ denote the degree of the polynomial and
\begin{align*}
&   \phi_1(g) = G_{22}(g), \\
&    \phi_2(g)=G_{22}(g)G_{11}(g)
-g^2\rho_0\left(C_{s1}^2k_{21}^2+C_{s2}^2k_{22}^2\right).
\end{align*}

Based on the characteristic polynomial derived through von Neumann analysis, 
we now proceed with the proof of Theorem \ref{thm:stability}, established in Section \ref{sec:semi_implicit}, which states the stability condition for the semi-implicit method. 
To do so, we first introduce
the following definitions and lemma \cite{strikwerda2004finite}. 

\begin{definition1}
A polynomial $ \phi $ is a von Neumann polynomial if and only if all its roots $ r_v $ satisfy the condition:
\begin{equation*}
    |r_v| \leqslant 1.
\end{equation*}
\end{definition1}

\begin{definition1}
A polynomial $ \phi $ is a simple von Neumann polynomial if and only if $ \phi $ is a von Neumann polynomial and its roots on the  unit circle are simple roots. 
\end{definition1}

\begin{definition1}
For a polynomial $ \phi $ of degree $ d $, its conjugate polynomial $ \phi^* $ is defined as:
\begin{equation*}
    \phi^*(g) = \sum_{l=0}^{d} \bar{a}_{d-l} g^l,
\end{equation*}
where $ \bar{a}_{d-l} $ denotes the conjugate of $ a_{d-l} $, moreover we define
\begin{equation*}
    \phi_{d-1}(g) = \frac{\phi_d^*(0)\phi_d(g) - \phi_d(0)\phi_d^*(g)}{g}.
\end{equation*}
\label{def:conj_poly}
\end{definition1}

\begin{lemma1}\label{lem:simple}
    $ \phi_d $ is a simple von Neumann polynomial if and only if \\
     (1) $ |\phi_d(0)| < |\phi_d^*(0)| $, and \\
     (2) $ \phi_{d-1} $ is a simple von Neumann polynomial.
\end{lemma1}

With these definitions and Lemma \ref{lem:simple} established, 
we now proceed to prove Theorem \ref{thm:stability}. 
\begin{proof}
It suffices to demonstrate that  the roots of  characteristic polynomial $\phi_3(g)$ (cf. equation \eqref{eq:charac_poly}) satisfy 
$|g|\leqslant 1$ under the CFL condition \eqref{eq:stab_cfl}. 
For the root of polynomial $\phi_1(g)$, denoted as $g_1$, we have the following inequality 
\begin{equation}
    |g_1|=\frac{|G_k|}{|1-C_{\mu1} k_{41}-C_{\mu2} k_{42}|}=\frac{|G_k|}{1+C_{\mu1} (1-\cos{\theta_1})+C_{\mu2} (1-\cos{\theta_2})}\leqslant |G_k|.
\end{equation}
The right side of this inequality can be shown to satisfy
\begin{align}
|G_k|\leqslant 1,
\label{eq:ineq_Gk}
\end{align}
provided that $C_m\in[0,1]$, where $C_m = C_{m1} + C_{m2}$. 

The proof will be complete once we establish that $\phi_2(g)$ is a simple von Neumann polynomial under the condition $C_m\in[0,1]$. 
Let $\phi_2(g)$ be expressed as 
\begin{equation*}
    \phi_2(g)=Ag^2+Bg+C,
\end{equation*}
where the coefficients are defined as:
\begin{equation*}
\begin{split}
A&= \left(1+ \rho_0 C_{s1}^2 C_{\mu3} k_{31}+\rho_0 C_{s2}^2 C_{\mu3} k_{32}\right)\left(1-C_{\mu1} k_{41}-C_{\mu2} k_{42}\right)-\rho_0 \left(C_{s1}^2k_{21}^2 + C_{s2}^2k_{22}^2\right) \\
&=(1+D)(1+E)
+\rho_0\left(C_{s1}^2\sin^2{\theta_1}+C_{s2}^2\sin^2{\theta_2}\right) > 1,\\
B&= G_k
\left(2+\rho_0 C_{s1}^2 C_{\mu3} k_{31}+\rho_0 C_{s2}^2 C_{\mu3} k_{32}-C_{\mu1} k_{41}-C_{\mu2} k_{42}\right)\\
&=G_k
(2+D+E),\\
C&= G_k^2,\\
D&= \rho_0 C_{s1}^2 C_{\mu3} (\cos{\theta_1}-1)^2+\rho_0 C_{s2}^2 C_{\mu3} (\cos{\theta_2}-1)^2, \\
E&= C_{\mu1} (1-\cos{\theta_1})+C_{\mu2} (1-\cos{\theta_2}).
\end{split}
\end{equation*}
The conjugate polynomial $ \phi_2^*(g) $ is given by (cf. definition \ref{def:conj_poly}):
\begin{equation}
    \phi^*_2(g) = \bar{C} g^2 + \bar{B} g + A.
\end{equation}
Now, we apply the first condition in Lemma \ref{lem:simple}, which states that 
 $\left|\phi_2(0)\right| < |\phi_2^*(0)|$. This implies  
\begin{equation}
 |C| < |A|.
\end{equation}
From equation \eqref{eq:ineq_Gk}, we know that 
 $|C| = |G_k^2| = |G_k|^2 \leqslant 1$. 
Additionally, since $A > 1$, the first condition of Lemma \ref{lem:simple} is satisfied. 

Next, we verify the second condition in Lemma \ref{lem:simple}, which requires: 
\begin{equation}
|AB - \bar{B}C| \leqslant A^2 - |C|^2. 
\end{equation}
Since $\bar{B}C = B|C|$, this condition becomes
\begin{equation}
|B| \leqslant A + |C|.
\end{equation}
Substituting the expressions for $A,B,C$, we obtain
\begin{align*}
    |G_k|&(2 + D+E)
    \leqslant (1 + D)  (1 + E) 
    + \rho_0 C_{s1}^2 \sin^2{\theta_1} + \rho_0 C_{s2}^2 \sin^2{\theta_2} + |G_k|^2.
\end{align*}
This inequality can be verified using the following relations: 
\begin{align*}
2|G_k| & \leqslant 1 + |G_k|^2, \\
    |G_k| (D+E) & \leqslant D+E,\\
    0 & \leqslant DE, \\
    0 & \leqslant \rho_0 C_{s1}^2 \sin^2{\theta_1} + \rho_0 C_{s2}^2 \sin^2{\theta_2}.
\end{align*}

Thus, by Lemma \ref{lem:simple}, we conclude that $\phi_2(g)$ is a simple von Neumann polynomial under the condition $C_m\in[0,1]$. 

\end{proof}

 \section{Iterative solver}
 \label{appendix:iter_solver}
 
 We propose an iterative procedure to solve the global system \eqref{eq:mass_semi_implicit}-\eqref{eq:rhie_chow}, following the framework of the SIMPLE method. The approach consists of three main steps: first, the momentum equation is solved using the pressure from the previous subiteration to obtain a predicted velocity; next, a linear system for the pressure correction is derived from the mass conservation equation; finally, the velocity field is corrected based on the computed pressure correction. The detailed procedure is outlined in Algorithm \ref{alg:semi_implicit_iter}.

Compared to the SIMPLE method, the proposed approach explicitly treats the nonlinear convective term, resulting in a linear system for the predicted velocity and another for the pressure correction that remain unchanged across subiterations and time steps. Consequently, these systems can be factorized in advance, significantly improving computational efficiency.
   
   \begin{algorithm}[H]                      
\caption{Iterative semi-implicit method}
\label{alg:semi_implicit_iter}

\begin{algorithmic}[1]
\State Initialize 
$(\mathbf{u}_{\rm f})_C^{n+1,0}, P_C^{n+1,0}$. 
\vspace{3pt}
\For {$k=0, \ldots, K$ }
    \State Compute the predicted velocity at cell center $(\mathbf{u}_{\rm f})_C^{n+1,k+\frac12}$: 
   \begin{equation*}
a_{C} (\mathbf{u}_{\rm f})_C^{n+1,k+\frac12} + \sum\limits_{F\in \text{nb}({C})} a_{F} (\mathbf{u}_{\rm f})_F^{n+1,k+\frac12} = -\nabla P_{C}^{n+1,k} + \mathbf{B}^n_C.
\end{equation*}
    
    \State Compute the predicted velocity at cell  boundary based on the Rhie et Chow interpolation (cf. equation \eqref{eq:rhie_chow}): 
 \begin{equation*}
 (\mathbf{u}_{\rm f})_b^{n+1,k+\frac12}=\frac12\left[(\mathbf{u}_{\rm f})_{C}^{n+1,k+\frac12}+(\mathbf{u}_{\rm f})_{C'}^{n+1,k+\frac12}\right]-(\mathbf{D}_{\rm f})_b\left[(\nabla P)_b^{n+1,k}-\overline{(\nabla P)_b^{n+1,k}}\right].
 \end{equation*}
 
 \State Obtain the pressure correction $P^{'}$ by solving the mass equation 
 \begin{equation*}
  \sum_{b(C)}\left[(\mathbf{u}_{\rm f})_{b}^{n+1,k+\frac12}
  +(\mathbf{u}_{\rm f})_b^{'}
  \right]\cdot\mathbf{S}_{\mathrm f,b}=0,
 \end{equation*}
\Statex \hspace{1.5em}  where the velocity correction $(\mathbf{u}_{\rm f})_b^{'}$ is related to the pressure correction $P_C^{'}$ by
$
 (\mathbf{u}_{\rm f})_b^{'}=-(\mathbf{D}_{\rm f})_b(\nabla P)_b^{'}.
$ 
 
\State  Determine the corrected velocity and pressure as
\begin{align*}
(\mathbf{u}_{\rm f})_{C}^{n+1,k+1} &= (\mathbf{u}_{\rm f})_{C}^{n+1,k+\frac12} -(\mathbf{D}_{\rm f})_C(\nabla P)_C^{'}, \\
P_{C}^{n+1,k+1} &= P_{C}^{n+1,k} + P_{C}^{'}. 
\end{align*}
 
    \State Terminate if the convergence criterion 
    $
    \frac{\Vert (\mathbf{u}_{\rm f})^{n+1,k+1} - (\mathbf{u}_{\rm f})^{n+1,k} \Vert}{\Vert \mathbf{u}_{\rm f}^{n+1,k} \Vert} < \texttt{tol}
    $
    is satisfied.
    \State Update the variables as follows:
    $
    (\mathbf{u}_{\rm f})_{C}^{n+1,k} \gets (\mathbf{u}_{\rm f})_{C}^{n+1,k+1}, \quad
    P_{C}^{n+1,k} \gets P_{C}^{n+1,k+1}.
    $
\EndFor
\end{algorithmic}
\end{algorithm}

\begin{remark*}
In the previous algorithm, we can set $K = 0$, meaning that no iterative procedure is involved. In this case, it can be straightforwardly derived that the approach is equivalent to solving
\begin{align*}
& \sum_{b(C)}(\mathbf{u}_{\rm f})_{b}^{n+1}\cdot\mathbf{S}_{\mathrm f,b}=0, \\
& a_{C} (\mathbf{u}_{\rm f})_C^{n+1} + \sum\limits_{F\in \text{nb}({C})} a_{F} (\mathbf{u}_{\rm f})_F^{n+1,\frac12} = -\nabla P_{C}^{n+1} + \mathbf{B}^n_C, \\
& (\mathbf{u}_{\rm f})_b^{n+1}=\frac12\left[(\mathbf{u}_{\rm f})_{C}^{n+1}+(\mathbf{u}_{\rm f})_{C'}^{n+1}\right]-(\mathbf{D}_{\rm f})_b\left[(\nabla P)_b^{n+1}-\overline{(\nabla P)_b^{n+1}}\right].
\end{align*}
Specifically, the mass conservation equation and the Rhie-Chow interpolation remain unchanged, while the momentum equation is solved using the predicted velocity for the neighboring cells of cell $C$. 

For certain scenarios, stable and reliable numerical results can still be obtained, provided that the CFL condition (cf. equation \eqref{eq:stab_cfl}) is satisfied. This approach offers efficiency advantages by reducing computational cost, as it eliminates the need for iterative corrections. In the numerical results, we explicitly indicate whether this semi-implicit method with $K=0$ is applied.

\end{remark*}

\bibliography{all_refs}

\end{document}